%% file: paper1.tex
\documentclass[final,5p]{elsarticle}

\usepackage{setspace}

\usepackage{natbib} 
\usepackage{amsmath,amsfonts,amssymb,amsbsy,amscd}
\usepackage{graphicx}
\usepackage{mathrsfs,array}
\usepackage{verbatim} 
\usepackage{subcaption}
\usepackage{color}

\usepackage[section,subsection,subsubsection]{placeins}

\newtheorem{thm}{Theorem}
\newtheorem{lem}[thm]{Lemma}
\newdefinition{rmk}{Remark}
\newproof{pf}{Proof}
\newproof{pot}{Proof of Theorem \ref{thm2}}

\newcommand{\Gammafs }{{\Gamma_{\!\!f\!s}}}

\newcommand{\II}{|\!|} 

\input{declarations}

\newcommand{\figopt}{figs3}
\renewcommand{\figopt}{figs_bw}

\newcommand{\vertiii}[1]{{\left\vert\kern-0.25ex\left\vert\kern-0.25ex\left\vert #1 
    \right\vert\kern-0.25ex\right\vert\kern-0.25ex\right\vert}}
\begin{document}

\begin{frontmatter}

%\title{Isogeometric analysis of water waves using linear potential theory} 
\title{Isogeometric analysis of linear free-surface potential flow}
%\title{Isogeometric analysis of potential flow with linearized free-surface}
%\title{Isogeometric potential flow with linearized free-surface problems}

\author{I. Akkerman}
\author{J.H.A. Meijer}
\author{M.F.P. ten Eikelder}
\address{ Mechanical, Maritime and Materials Engineering Department  \\ 
 Delft University of Technology}
\date{\today}

\begin{abstract}
This paper presents a novel variational formulation to simulate linear free-surface flow. The variational formulation is 
suitable for higher-order finite elements and higher-order and higher-continuity shape functions as employed in Isogeometric Analysis (IGA).

%The method employs a monolithic approach to solve the potential and the water height.
The novel formulation combines the interior and free-surface problems in one monolithic formulation.
This leads to exact energy conservation and superior performance in terms of accuracy when compared to a traditional segregated  formulation.
This is confirmed by the numerical computation of traveling waves in a periodic domain and a three-dimensional sloshing problem.
The isogeometric approach shows significant improved performance compared to traditional finite elements.
Even on very coarse quadratic NURBS meshes the dispersion error is virtually absent.
\end{abstract}


\end{frontmatter}

%=============================================
\section{Introduction}
\label{sec:intro}
%=============================================
In maritime and offshore applications potential flow is often 
a very useful model for prescribing relevant design scenarios.
This is particularly the case when waves, which are dominated by 
inertia effects, are of main concern.

Numerical methods for predicting potential flow are often based on 
boundary integral formulations, referred to as boundary element methods (BEM) or panel methods. 
Examples of successful application of the BEM are (i) wave drift forces \cite{Pinkster79},
(ii) ship wave resistance \cite{Dawson77} and (iii) ship motion and maneuvering \cite{STF70}. 
These references indicate the starting point of the development of BEM and numerous works have been published 
to improve these methods.

All these methods exploit the Greens functions and only require  the boundary to be  discretized. 
This leads to very low degree-of-freedom counts, however, the resulting matrices are dense and expensive to assemble.
Moreover, straightforward implementations have a very unfavorable scaling in terms of computational time versus degree-of-freedom. 

One way to circumvent this scaling issue is the use of the fast multipole method \cite{CRW93}. 
An alternative approach, which is explored in this paper, is to forgo the boundary integral approach.
Instead variational methods discretizing the entire volume are used, this encompasses the classical finite element method (FEM) 
as well as the novel Isogeometric Analysis (IGA) approach \cite{HuCoBa04,CoHuBa09}.
IGA uses splines as shape functions as opposed to the more standard polynomials in FEM.
These spline functions are also used in CAD. One of the main goals of IGA is to create a seamless integration of analysis  
into design processes by using the same geometry description. 
It turns out that the spline shape functions have several other benefits, 
in terms of efficiency and accuracy \cite{BACHH07,CoHuRe07,EBBH09}.

Linear finite elements are used for the simulation of water waves by for instance
Wu and Eatock-Taylor  \cite{WuTaylor94},
Kim and Bai {\it et al} \cite{kim1999finite},
and  Westhuis \cite{phdWesthuis}, just to name a few.
In the first two publications it is the starting point of a large body of work.
An important issue in these methods is the coupling of the interior Laplace problem with the free-surface evolution, in particular the horizontal velocity at the free-surface. 
The mentioned bodies of work deal with this issue in different manners.
In Wu and Eatock-Taylor \cite{WuTaylor94} the recovered velocity is a projection on a new finite element space.
This is referred to as Galerkin projection and results in a mixed formulation.
Kim and Bai \cite{kim1999finite} employ a  similar potential-velocity mixed formulation but employ the Hamiltonian formalism.  
In Westhuis \cite{phdWesthuis} the vertical velocity at the free-surface is reconstructed using finite differences. 
Zienkiewicz and Zhu \cite{ZZ92} have used other reconstruction methods leading to adequate results.

The previous works do allow for nonlinear waves.
In this work we take a step back and limit ourselves to linear waves. 
However, we perform a thorough investigation of the numerical behavior of the methods considered
with a particular focus on the coupling between the interior and free-surface problem.
We investigate several weak formulations that deal with the coupling in different manners. 
Moreover, we establish the well-posedness of these formulations.
In the discrete case, we derive accuracy estimates and analyze the energy behavior. 
Using an appropriate time integrator for some weak formulations provides exact energy conservation.
The spatial discretization employs higher order elements.

The outline of the paper is as follows.
In section \ref{sec:strong} the strong formulation of the problem is introduced.
In section \ref{sec:weak} several corresponding weak formulations are presented 
and their energy conservation properties are analyzed.
Energy conservation is only guaranteed when an appropriate time integrator is used. In section  
\ref{sec:time_int}  midpoint time integration is shown to achieve this.
After the temporal discretisation, we will focus on the spatial discretisation in section \ref{sec:spatial}.
In this section particular attention is given to Isogeometric analysis.
In section \ref{sec:numres} the energy conservation  and dispersion properties of the methods 
are investigated.
This reveals the superiority of the monolithic approach over the segregated approach. 
In section \ref{sec:coerc} and \ref{sec:VV}, the convergence of monolithic formulation is 
analyzed and verified, respectively. 
Finally, the paper ends with a 3D show case problem in section \ref{sec:show}, and conclusions in section \ref{sec:conc}.

%=============================================
\section{Strong forms of the free-surface problem}
\label{sec:strong}
%=============================================

In this section we introduce three different sets of governing equations describing potential waves.
We start off with the nonlinear equations, which are then linearized and condensed to arrive at two alternative formulations.

%=============================================
\subsection{Non-linear strong form}
\label{ssec:nonlin}
%=============================================
The governing equations read:
\begin{subequations}\label{eq:strong_nonlin}
\begin{align}
\Delta \phi =&~0 && \text{in~} \Omega, \label{eq: LP}\\
\phi_  +\onehalf\nabla \phi \cdot \nabla\phi+  g \eta =& ~0 && \text{on~} \Gammafs, \label{eq: dyn} \\
\eta_t+ \widetilde{\nabla} \phi \widetilde{\nabla} \eta - \phi_{z}  =&~0&& \text{on~} \Gammafs, \label{eq: kin}
\end{align}
\end{subequations}
where $\phi:\Omega \rightarrow \mathbb{R}$ is the potential,  $\eta:\Gammafs \rightarrow \mathbb{R}$ is the water height, and
$\widetilde{\nabla}$ is the gradient restricted to exclude the vertical direction. Here $\Gamma_{fs}$ denotes the free-surface.

The Laplace problem (\ref{eq: LP}) is a consequence of the conservation of mass
%\begin{align}
%\frac{\partial \rho}{\partial t}  + \nabla \cdot \rho \bu =&0 && \text{in~} \Omega
%\end{align}
for a constant density fluid, 
\begin{align}
\nabla \cdot \bu =&~0 \quad \text{in~} \Omega
\end{align}
in conjunction with the assumption of an inviscid and irrotational fluid, for which we can write
\begin{align}
\bu =&~\nabla \phi.
\end{align}

The dynamic boundary condition (\ref{eq: dyn}) is a condition on the pressure, 
whereas the kinematic condition (\ref{eq: kin}) ensures that the surface moves with the water.

 The problem needs to be augmented
with boundary conditions on the remaining boundaries.
In this work we will assume either periodic boundary conditions or no-penetration boundary conditions.
In the first case no boundary term is present, as such no boundary condition needs to be enforced.
In case of the no-penetration boundary conditions, we get a Neumann boundary condition
\begin{align}
\bu \cdot \bn = \bn\cdot \nabla \phi =&~ 0 && \text{on~} \Gamma/\Gammafs.
\end{align}
The problem is completed by specifying appropriate initial conditions for $\phi$ and $\eta$ at the free-surface.
Note that in the fully non-linear problem, the location of the free-surface, 
denoted with $\Gammafs$, is part of the solution as it is determined by $z=\eta$.

%=============================================
\subsection{Linear strong form}
\label{ssec:linear}
%=============================================
The nonlinear equations (\ref{eq:strong_nonlin}) can be simplified by assuming small disturbances.
By neglecting the quadratic terms we arrive at
\begin{subequations}\label{eq:strong_lin}
\begin{align}
\Delta \phi =&~0 && \text{in~} \Omega,\\
\phi_t  + g \eta =& ~0 && \text{on~} \Gammafs,\label{eq:strong_dyn}\\
\phi_z =& ~\eta_t && \text{on~} \Gammafs,\label{eq:strong_kin}\\
\bn\cdot \nabla \phi =& ~0 && \text{on~} \Gamma/\Gammafs,
\end{align}
\end{subequations}
where $\Gammafs$ is now assumed frozen on the undisturbed location.

By combining the dynamic and kinematic boundary condition, (\ref{eq:strong_dyn}) and (\ref{eq:strong_kin}), we can eliminate the water height $\eta$ from the problem.
The problem reduces to:
\begin{subequations}\label{eq:strong_comb}
\begin{align}
\Delta \phi =&~0&& \text{in~} \Omega,\label{eq:strong_comb_a}\\
\phi_{tt}  + g \phi_{z}  =& ~0 && \text{on~} \Gammafs,\label{eq:strong_comb_b}\\
\bn\cdot \nabla \phi =& ~0 && \text{on~} \Gamma/\Gammafs,\label{eq:strong_comb_c}
\end{align}
\end{subequations}
where the problem is now second order in time. This means that the initial condition for $\eta$ at the free-surface
is replaced by an initial condition for $\phi_t$ at the free-surface.

%=============================================
\section{Weak forms of the linear free-surface problem}
\label{sec:weak}
%============================================
In this section we present several weak formulations of the linear wave problem.
We start with a weak form of (\ref{eq:strong_comb}), and subsequently propose several formulations for problem (\ref{eq:strong_lin}).
These later formulations of mixed character are more amiable in the nonlinear case or in situations with currents present.
Furthermore, we analyze the  energy conservation of the solution for each of the formulations.
An analysis of the existence and accuracy  of the solution for each of the formulations, is postponed to section \ref{sec:coerc}.

We introduce the notation
\begin{subequations}
\begin{alignat}{1}
(w,v)_\Omega  =& \int_\Omega \!\! wv \; d\Omega, \\
(w,v)_\Gammafs   =& \int_\Gammafs \!\! wv \; d{\Gamma} ,
\end{alignat}
\end{subequations}
to denote the innerproduct over the entire domain and free-surface, respectively.
The corresponding norms are denoted with
\begin{subequations}
\begin{alignat}{1}
\|w\|^2_\Omega   =~&(w,w)_\Omega,\\
\|w\|^2_\Gammafs =~&(w,w)_\Gammafs.
\end{alignat}
\end{subequations}
%=============================================
\subsection{Weak formulation  of the reduced problem}\label{sec:eng1}
%=============================================
%The weak form of (\ref{eq:strong_comb}) can be obtained by 
%multiplying (\ref{eq:strong_comb_a}) with an arbitrary function ($w$). 
%The order of the required derivative can be reduced by using the divergence theorem, the resulting 
%boundary terms can be simplified by using (\ref{eq:strong_comb_b}) and  (\ref{eq:strong_comb_c}).

A weak form of (\ref{eq:strong_comb}) follows when multiplying (\ref{eq:strong_comb_a}) with an arbitrary function and integrating over the domain. 
The order of the required derivatives can be reduced by using Green's identities and the resulting 
boundary terms can be simplified by using (\ref{eq:strong_comb_b}) and  (\ref{eq:strong_comb_c}).
Let $\WW=H^1(\Omega)$ denote the function-space.
The variational formulation of the reduced problem reads:\\
\vspace{1.2mm}\\
\textit{Find $\phi \in \WW$ such that for all $ w \in \WW$:}
\begin{align} \label{eq:form1}
B_r(w,\phi) = 0
\end{align}
where 
\begin{align}
B_r(w,\phi)  :=
(\nabla w,\nabla \phi)_\Omega +\frac{1}{g} (w, \phi_{tt} )_{\Gammafs}.
\end{align}

An energy conservation statement for the reduced form can be derived by choosing $w=\phi_t$ in (\ref{eq:form1}). This selection gives
\begin{align}
B_r(\phi_t,\phi)   =& (\nabla \phi_t,\nabla \phi)  +\frac{1}{g} (\phi_t, \phi_{tt} )_{\Gammafs}  \nonumber \\
=&\dfrac{\rm d}{{\rm d}t} \onehalf \| \nabla \phi \|^2  + \dfrac{\rm d}{{\rm d}t} \frac{1}{2 g} \|\phi_t\|^2_{\Gammafs} =0.
\end{align}
Realizing that the definition for kinetic and potential are 
\begin{subequations}
\begin{alignat}{1}
E_{\rm kin} =& \onehalf \| \bu \|^2  =  \onehalf \| \nabla \phi \|^2,\\
E_{\rm pot} =& \onehalf g \| \eta \|^2 =  \frac{1}{2g} \|\phi_t\|^2_{\Gammafs},
\end{alignat}
\end{subequations}
we arrive at the following statement
\begin{align}
 \dfrac{\rm d}{{\rm d}t} E_{\rm kin} +  \dfrac{\rm d}{{\rm d}t} E_{\rm pot} = 0.
\end{align}
Consequently, the reduced weak form is exactly energy conservative assuming appropriate time integration.

%=============================================
\subsection{Weak formulations of the segregated problem}
%=============================================

Here we present a segregated weak formulation serving as a reference method based on the work of Wu and Taylor \cite{WuTaylor94} and Wu et al. \cite{WuMaTaylor98},
Westhuis \cite{phdWesthuis}, Kyoung et al. \cite{KHKB2005}, Bai et al. \cite{BCCK05} and Kim et al. \cite{KKEB2006}.
It is based on the strong form (\ref{eq:strong_lin}) and decouples the interior and surface parts of the problem.

The interior problem in strong form reads:
\begin{subequations}
\begin{alignat}{3}
\Delta \phi =&0 &\qquad& \text{in~} \Omega,\\
\phi  =& \hat{\phi} &\qquad& \text{on~} \Gammafs,\\
n\cdot \nabla \phi =& 0 &\qquad& \text{on~} \Gamma/\Gammafs,
\end{alignat}
\end{subequations}
where $\hat{\phi}$ is input from the free--surface problem:
\begin{subequations}
\begin{alignat}{3}
\hat{\phi}_t  + g \eta =& 0 &\qquad& \text{on~} \Gammafs,\\
 \eta_t=&\phi_z&\qquad& \text{on~} \Gammafs.
\end{alignat}
\end{subequations}
Here is $\phi_z$ is given by the interior problem. %This interior problem will give $\hat{\phi}$ for the next time step.
As such the two problems are artificially decoupled.

Let $\WW_0$ and $\WW_{\hat{\phi}}$ denote the subspaces of $H^1(\Omega)$ satisfying homogeneous Dirichlet boundary condition, $\phi = 0$ on $\Gammafs$, and the homogeneous Dirichlet boundary condition $\phi = \hat{\phi}$ on $\Gammafs$, respectively. Furthermore let $\WW_\Gammafs$ denote the trace space of $\WW$.
% and set $\VV_\Gammafs = \WW_\Gammafs \times \WW_\Gammafs$.
%We introduce the notation $\mathbf{W}=(w,v)$ and $\bvphi=(\phi, \eta)$.
The associated weak formulation for the interior problem reads:\vspace{1.2mm}\\
\textit{Find $\phi \in \WW_{\hat{\phi}}$ such that for all $ w \in \WW_0$:}
\begin{align}\label{eq:form2a}
B_{int}(w,\phi) = 0
\end{align}
where 
\begin{align}
B_{int}(w,\phi)  := (\nabla w,\nabla \phi). 
\end{align}
For the free-surface the weak form reads:\vspace{1.2mm}\\
\textit{Find $(\phi,\eta) \in \WW_{\Gamma}\times\WW_{\Gamma}$ such that for all $ (w,v) \in \WW_{\Gamma}\times\WW_{\Gamma}$:}
\begin{align}\label{eq:form2b}
B_{fs}(\{w,v\};\{\phi,\eta\}) =  \frac{g^2}{\alpha^2} (v ,\phi_z)_{\Gammafs},
\end{align}
where
\begin{align}
B_{fs}(\{w,v\};\{\phi,\eta\}) := 
(w,\phi_t  + g \eta)_{\Gammafs}
+\frac{g^2}{\alpha^2} (v , \eta_t)_{\Gammafs}.
\end{align}
Here $\alpha$ is a  parameter, that eventually will be chosen to depend on the time integrator. % that has the unit $s^{-1}$

\subsubsection{Energy conservation of the segregated formulation}

To arrive and an energy statement we would like to select $w=\phi_t$ in both the interior problem (\ref{eq:form2a})
 and free-surface problem (\ref{eq:form2b}). However,  due to the Dirichlet boundary condition on $\Gammafs$ 
this is not allowed for the interior problem. To remedy this, the weak form (\ref{eq:form2a}) is written in an equivalent Lagrange multiplier formulation. 
from which yields exactly the same solution if the Lagrange multiplier space is constructed appropriately.
This Lagrange multiplier formulation reads as follows:

\textit{Find $(\phi,\lambda) \in \WW\times \WW_{\Gamma}$ such that for all $ (w,q)\in \WW\times \WW_{\Gamma}$:}
\begin{align}\label{eq:form2a_LM}
(\nabla w,\nabla \phi)_\Omega+ (w, \lambda)_{\Gammafs} + (q, \phi)_{\Gammafs} = (q, \hat{\phi})_{\Gammafs}.
\end{align} 
The new formulation allows the selection of $w=\phi_t$ and additionally we set $q=0$.
This yields:
\begin{align}
(\nabla \phi_t,\nabla \phi) + (\phi_t, \lambda) = 0,
\end{align}
which can be rewritten as:
\begin{align}\label{eq: segregated E kin}
\dfrac{\rm d}{{\rm d}t} E_{\rm kin} = - (\phi_t, \lambda). 
\end{align}

For the free-surface  problem we select $w=\phi_z$ and $v =\eta \alpha^2/g$ which results in:
\begin{align}
0=& (\phi_z,\phi_t  + g \eta)_{\Gammafs}
+g ( \eta , \eta_t - \phi_z)_{\Gammafs} \nonumber\\
=&(\phi_z,\phi_t )_{\Gammafs}
+g ( \eta , \eta_t )_{\Gammafs}.
\end{align}
This can be rewritten as:
\begin{align}
\dfrac{\rm d}{{\rm d}t} E_{\rm pot} = - (\phi_z,\phi_t )_{\Gammafs}.
\end{align}

Combining the interior and boundary problem we arrive at the following energy statement:
\begin{align}
\dfrac{\rm d}{{\rm d}t} E_{\rm kin}  +\dfrac{\rm d}{{\rm d}t} E_{\rm pot}  = - (\phi_t, \lambda +\phi_z) 
\end{align}
which results in conservation of energy if $\lambda = - \phi_z $.
Unfortunately this relation only holds on sufficiently smooth meshes, where the solution is converged.
This means conservation of energy can not be guaranteed.

In the following two subsections two alternative approaches to remedy this energy error are discussed.
The first approach is to reconstruct the Lagrange multiplier and use this as a forcing in the boundary formulation.
The other approach is to solve the boundary and interior problem in one monolithic formulation.
 
%=============================================
\subsection{Segregated formulation, with LM reconstruction}
%=============================================
The Lagrange multiplier from equation (\ref{eq:form2a_LM}) can be reconstructed 
by selecting $q=0$, which results in 
\begin{align}\label{eq:form2a_LM_recon}
(w, \lambda)_{\Gammafs}  = - (\nabla w,\nabla \phi)_\Omega 
\end{align} 
for all $ w\in \WW$.

Energy conservation can be recovered if this Lagrange multiplier is directly used in 
the weak form for the free-surface problem given in (\ref{eq:form2b}).
The right-hand side would need to be modified to:
\begin{align}
  \frac{g^2}{\alpha^2} (v ,\phi_z)_{\Gammafs}=
  - \frac{g^2}{\alpha^2} (v ,\lambda)_{\Gammafs}=
 \frac{g^2}{\alpha^2} (\nabla v ,\nabla \phi),
\end{align}
where in the last integral the functions $v$ need to be arbitrarily extended into the domain.

For the modified free-surface problem we again select $w=\phi_z$ and $v =\eta \alpha^2/g$ which results in:
\begin{align}
\dfrac{\rm d}{{\rm d}t} E_{\rm pot} 
  = (\lambda,\phi_t)_{\Gammafs}.
\end{align}
This can be combined with already obtained kinetic energy statement (\ref{eq: segregated E kin}) to yield:
\begin{align}
\dfrac{\rm d}{{\rm d}t} E_{\rm kin} +\dfrac{\rm d}{{\rm d}t} E_{\rm pot} =0,
\end{align}
and as such recovering the energy conservation.

%Note that the in the work of Westhuis \cite{phdWesthuis}  a higher-order reconstruction of $\phi_z$ is used.
%This reconstruction probably solves some of the problems highlighted here.
\begin{rmk}
Note that in the work of Westhuis \cite{phdWesthuis}  a higher-order reconstruction of $\phi_{z} $ is used.
This reconstruction does not provide exact energy conservation.
\end{rmk}

%=============================================
\subsection{Monolithic weak formulation}
%=============================================
Here we present a monolithic formulation that bypasses the need of the Lagrange multiplier construction. 
When the divergence theorem on the interior problem creates boundary terms 
the kinematic boundary condition can be substituted in, this eliminates the problematic derivative. 
The dynamic boundary condition is added to the weak form in a way that guarantees coercivity.%and obtain an optimal error estimate.

%Let the test function pair be $\bW:=(w,v) \in \VV$ and the trial function pair denote $\bvphi:=(\phi,\eta) \in \VV$.
We propose the weak formulation:\\
\textit{Find $(\phi,\eta) \in \WW\times\WW$ such that for all $ (w,v) \in \WW\times\WW$,}
\begin{align}
B(\{w,v\};\{\phi,\eta\}) = 0,
\end{align}
with 
\begin{align} \label{eq:form_mono}
B(\{w,v\};\{\phi,\eta\}) =& 
(\nabla w,\nabla \phi)_\Omega - (w, \eta_t)_{\Gammafs}  \nonumber \\ &
+\onehalf (v + \textstyle{\frac{\alpha}{g}} w,\phi_t  + g \eta)_{\Gammafs}.
\end{align}
Again $\alpha$ is a parameter that eventually will be chosen based on the time integrator. % that has the unit $s^{-1}$

\subsubsection{Conservation of energy  of the monolithic form}\label{sec:eng4}
To establish an energy statement we select
$w= \phi_t$  and $v =2\eta_t -  \frac{\alpha}{g} \phi_t$
and substitute this in (\ref{eq:form_mono})
which gives:
\begin{align}
B(\{\phi_t,\eta_t\};\{\phi,\eta\}) =& 
(\nabla \phi_t,\nabla \phi)  - (\phi_t, \eta_t)_{\Gammafs}  \nonumber \\
&+\onehalf (2\eta_t -  \textstyle{\frac{\alpha}{g}} \phi_t +  \textstyle{\frac{\alpha}{g}} \phi_t,\phi_t  + g \eta)_{\Gammafs} \nonumber  \\
%=& (\nabla \phi_t,\nabla \phi)  - (\phi_t, \eta_t)_{\Gammafs} 
%+   (\eta_t , \phi_t  + g \eta)_{\Gammafs}   \nonumber \\
=& 
(\nabla \phi_t,\nabla \phi)  
+ g  (\eta_t , \eta)_{\Gammafs}.
\end{align}
From this it follows directly that
%\begin{align}
%E_{\rm kin} =& \onehalf \| \bu \|^2  =  \onehalf \| \nabla \phi \|^2 \\
%E_{\rm pot} =& \onehalf g \| \eta \|^2 
%\end{align}
\begin{align}
\dfrac{\rm d}{{\rm d}t}E_{\rm kin} + \dfrac{\rm d}{{\rm d}t}E_{\rm pot} =0,
\end{align}
which indicates that the total energy is conserved.

%=============================================
\section{Time integration}\label{sec:time_int}
%=============================================
In sections 
\ref{sec:eng1} and
\ref{sec:eng4}
it was proven for the respective weak forms that the total energy is conserved. This is stated as
\begin{align}
\dfrac{\rm d}{{\rm d}t}E_{\rm kin} + \dfrac{\rm d}{{\rm d}t}E_{\rm pot} =0.
\end{align}
When defining the energies as
\begin{subequations}
\begin{alignat}{1}
E^n_{\rm kin} =&  \onehalf (\phi^n, \phi^n),  \\
E^n_{\rm pot} =& \onehalf g (\eta^n, \eta^n),
\end{alignat}
\end{subequations}
and assuming the time-integrator has the correct behavior the time continous statement can be translated to a time discrete conservation statement, namely
\begin{align}
 E^{n+1}_{\rm kin} + E^{n+1}_{\rm pot} =& E^{n}_{\rm kin} + E^n_{\rm pot}.
\end{align}
This  translation holds when the time integrator satisfies the following relations,
\begin{subequations}\label{eq:timeint_req}
\begin{alignat}{1}
\frac{(\phi^{n+1}, \phi^{n+1})- (\phi^{n}, \phi^{n})}{2 \Delta t} =&  (\tilde{\phi}, \tilde{\phi}_t),
 \label{eq:kin_req}\\ 
g\frac{(\eta^{n+1}, \eta^{n+1})-  (\eta^{n}, \eta^{n})}{2 \Delta t}  =& g (\tilde{\eta}, \tilde{\eta}_t).  \label{eq:pot_req}\
\end{alignat}
\end{subequations}
Here the tilde denotes the value used by the time integrator to evaluate  the weak formulation.

In this paper midpoint time integration is adopted for both the first and second-order problems.
In both cases the required translation is valid, as will be shown, while having second-order convergence.

Note, that as shown in \cite{EiAk17i} generalized midpoint time integration almost satisfies (\ref{eq:timeint_req}). 
In addition to (\ref{eq:timeint_req}) it also features a diffusion that scales with $\mathcal{O}(\Delta t)$.

%a one parameter family of time integration routines is derived that guarantees diffusion.
%The one parameter family is generated using a specific selection of parameters in the generalized-$alpha$ time integration routine.
%Note that this family of integrators  

%=============================================
\subsection{Midpoint time integration for a first-order problem}
%=============================================
For problems with first-order time derivatives the midpoint time-integration is determined by the following relations.
\begin{itemize}
\item The ordinary differential equation:
\begin{align}
f(\phi_t^{n+1/2 },\phi^{n+1/2 } ) = 0.
\end{align}
\item  The interpolation relation:
\begin{align}\label{eq:1st-int}
\phi^{n+1/2 } =& \onehalf \left ( \phi^{n+1} + \phi^{n} \right ) .
\end{align}
\item  The  kinematic relation:
\begin{align}\label{eq:1st-kin}
{\phi}^{n+1 }   = &  {\phi}_{n} + \Delta t {\phi_t}^{n+1/2 }.
\end{align}
\end{itemize}
This results in three relations for three unknowns that are in principle solvable.

In \ref{sec:coerc_mono} and \ref{app:coerc_seg}  it will be shown that using
\begin{align}
\alpha = \frac{\partial \tilde{\phi}_{t}}{\partial \tilde{\phi}}
\end{align}
results in favourable properties of the formulations.
For the midpoint time-integration this results in:
\begin{align}
\alpha = \frac{\partial \phi_{t}^{n+1/2 }}{\partial \phi^{n+1/2 }}
%  = \frac{\partial \phi_{t}^{n+\onehalf }}{\partial \phi^{n+1}}
 %  \frac{\partial \phi^{n+1}}{\partial \phi^{n+\onehalf }}
 %= \left ( \frac{\partial \phi^{n+1}}{\partial \phi_{t}^{n+\onehalf }} \right )^{-1}
% \left (  \frac{\partial \phi^{n+\onehalf}}{\partial \phi^{n+1}}\right )^{-1}
  = \left ( \frac{\partial \phi^{n+1}}{\partial \phi_{t}^{n+1/2 }}
 \frac{\partial \phi^{n+1/2}}{\partial \phi^{n+1}}\right )^{-1}
 % =  \left (  \Delta t  \right )^{-1} 2 
   = \frac{2}{\Delta t}.
 %=  \frac{\partial \phi^{n+\onehalf }}{\partial \phi_{tt} ^{n+\onehalf }} \right )^{-1}
 % = \left ( \frac{\partial \phi^{n+\onehalf }}{\partial \phi_{tt} ^{n+\onehalf }} \right )^{-1}
\end{align}

\subsubsection{Kinetic and potential  energy behavior}
To see whether the kinetic energy requirement (\ref{eq:kin_req}) holds,
we rewrite  the kinematic relation (\ref{eq:1st-kin}) as:
\begin{align}
   {\phi_t}^{n+1/2 }    = &    \frac{1}{\Delta t } \left ({\phi}^{n+1}   - {\phi}^{n}    \right ).
\end{align}
Using this relation we find that:
\begin{align}
(   \phi_t^{n+1/2 },\phi^{n+1/2 }  )   = &  \left (  \frac{1}{\Delta t } \left ({\phi}^{n+1}   - {\phi}^{n}    \right ), \onehalf (  {\phi}^{n} +  {\phi}^{n+1}   )\right )  \nonumber \\
%=& \frac{1}{2\Delta t }\left (  {\phi}^{n+1}   - {\phi}^{n},  {\phi}^{n} +  {\phi}^{n+1}   \right )\nonumber \\
=& \frac{1}{2\Delta t }\left (  \left (  {\phi}^{n+1}  ,    {\phi}^{n+1}   \right ) -\left (  {\phi}^{n}  ,    {\phi}^{n}   \right ) \right ).
\end{align}
This demonstrates that the kinetic energy requirement is satisfied.
The potential energy relation goes analogously. 

%=============================================
\subsection{Midpoint time-integration for a 2nd order problem}
%=============================================
For problems with 2nd order time derivatives the midpoint time-integration is determined by the following relations.
\begin{itemize}
\item The ordinary differential equation
\begin{align}\label{eq:2nd-pde}
f(\phi_{tt} ^{n+1/2 } ,\phi_t^{n+1/2 },\phi^{n+1/2 } ) = 0.
\end{align}
\item Two interpolation relations,
\begin{subequations}\label{eq:2nd-int}
\begin{alignat}{1}
   {\phi}^{n+1/2 }   = & \onehalf (  {\phi}^{n} +  {\phi}^{n+1}   ),\\
     {\phi_t}^{n+1/2 }    = & \onehalf (  {\phi_t}^{n} +  {\phi_t}^{n+1}   ).
\end{alignat}
\end{subequations}
\item Two kinematic relations,
\begin{subequations}\label{eq:2nd-kin}
\begin{alignat}{1}
{\phi}^{n+1 }   = &  {\phi}^{n} + \Delta t {\phi_t}^{n+1/2 },\label{eq:2nd-kin1} \\
  {\phi_t}^{n+1 }    = &    {\phi_t}^{n} +  \Delta t {\phi_{tt} }^{n+1/2 }. \label{eq:2nd-kin2}
\end{alignat}
\end{subequations}
\end{itemize}
This results in five relations for five unknowns that are in principle solvable.

\subsubsection{Kinetic energy behavior}
Again, to see whether the kinetic energy requirement  (\ref{eq:kin_req}) holds
we rewrite  the kinematic relation (\ref{eq:2nd-kin1}) as
\begin{align}
   {\phi_t}^{n+1/2 }    = &    \frac{1}{\Delta t } \left ({\phi}^{n+1}   - {\phi}^{n}    \right ).
\end{align}
Using this relation we again find that
\begin{align}
(   \phi_t^{n+1/2 },\phi^{n+1/2 }  )   = &  \left (  \frac{1}{\Delta t } \left ({\phi}^{n+1}   - {\phi}^{n}    \right ), \onehalf (  {\phi}^{n} +  {\phi}^{n+1}   )\right )  \nonumber \\
%=& \frac{1}{2\Delta t }\left (  {\phi}^{n+1}   - {\phi}^{n},  {\phi}^{n} +  {\phi}^{n+1}   \right )\nonumber \\
=& \frac{1}{2\Delta t }\left (  \left (  {\phi}^{n+1}  ,    {\phi}^{n+1}   \right ) -\left (  {\phi}^{n}  ,    {\phi}^{n}   \right ) \right ),
\end{align}
demonstrating that the kinetic energy requirement is satisfied.

\subsubsection{Potential energy behavior}
For the potential energy we need to make the following substitution,
%\begin{subequations}
%\begin{alignat}{1}
%\eta =& - \frac{1}{g}   \phi_t \\
%\eta_t =& - \frac{1}{g}  \phi_{tt}   %= -\frac{1}{\Delta t \;g } \left ({\phi_t}^{n+1}   - {\phi_t}^{n}    \right ).
%\end{alignat}
%\end{subequations}
\begin{align}
\eta =& - \frac{1}{g}   \phi_t 
\end{align}

which is justified by the dynamics boundary condition (\ref{eq:strong_kin}). % and  the kinematic relation (\ref{eq:2nd-kin2}).
The potential energy requirement  (\ref{eq:pot_req}) holds, namely,
\begin{align}
g (\eta^{n+1/2 }, \eta_t^{n+1/2 })
% = &\frac{1}{g} (   \phi_{tt} ^{n+1/2 },\phi_t^{n+1/2 }  )    \nonumber \\
= &   \frac{1}{g} \left (    \frac{1}{\Delta t } \left ({\phi_t}^{n+1}   - {\phi_t}^{n}    \right ), \onehalf (  {\phi_t}^{n} +  {\phi_t}^{n+1}   )\right )  \nonumber \\
%= &  \frac{1}{2g\Delta t }\left ( {\phi_t}^{n+1}   - {\phi_t}^{n}    ,   {\phi_t}^{n} +  {\phi_t}^{n+1}   \right )  \nonumber \\
= &  \frac{1}{2g\Delta t }\left ( \left ( {\phi_t}^{n+1} ,  {\phi_t}^{n+1}   \right ) -\left ( {\phi_t}^{n} ,  {\phi_t}^{n}   \right )   \right ) \nonumber \\
= &  \frac{g}{2\Delta t }\left ( \left ( \eta^{n+1} ,  \eta^{n+1}   \right ) -\left ( \eta^{n} ,  \eta^{n}   \right )   \right ).
\end{align}

%=============================================
\section{Spatial discretization}\label{sec:spatial}
%=============================================
In this section the spatial discretisation is discussed.
For the spatial discretisation both finite elements (FE) and NURBS based Isogeometric analysis (IGA) are used.
IGA can be seen as an extension of FE. Both methods will be explained separately. But first the commonalities are discussed.
Both approaches approximate the unknown exact solution, denoted as $\phi(\bx)$, by a weighted sum of known shape functions $N_b(\bx)$:
\begin{align}
\phi(\bx) \approx \phi^h(\bx) = \sum_{ b=1}^{n_{dof}}  \phi_b N_b(\bx).
\end{align}
Here  $\phi_b$ are the unknown parameters that need to be determined.
To convert the weak form in a set of equations we select:
\begin{align}
w(\bx)  =  N_a(\bx) \qquad ~a = 1,2,\dots, n_{dof}.
\end{align}
Using these approximations the terms in the weak form can be determined. For instance, 
a term in the weak form such as $(w ,\phi) $ results in:
\begin{align}
( \nabla N_a, \nabla \phi^h) 
=& \int   \nabla N_a \cdot 
%\qquad \qquad \nonumber \\& \qquad \qquad  
\left ( \sum_{b=1}^n \phi_b   \nabla N_b \right ) ~d\Omega
\nonumber \\
=&  \sum_{b =1}^n  \int   \nabla N_a   \cdot    \nabla N_b  ~d\Omega ~ \phi_b%\nonumber \\
= \sum_{b =1}^n K_{ab} ~ \phi_b,
%\nonumber \\
%=& \sum_{a=1}^n  \sum_{b =1}^n w_a \sum_{e =1}^{n_{el}} \int_{\Omega_e}  N_a(\bx)  \cdot   N_b(\bx) ~d\Omega ~ \phi_b
\end{align}
where we dropped the argument $\bx$ to simplify the notation.
The matrix is defined as:
\begin{align}
K_{ab} := \int   \nabla N_a  \cdot   \nabla N_b  ~d\Omega=  \sum_{e =1}^{n_{el}} \int_{\Omega_e}   \nabla  N_a   \cdot   \nabla N_b~d\Omega,
\end{align}
Here $\Omega_e$ denotes the element domain,of which the union covers the entire domain, viz. $\Omega = \cup_e \Omega_e$.
The element integral is approximated using Gauss quadrature. 
Both the domain of the integral and gradients incorporate the mapping between the physical and reference domain, following the isoparametric paradigm.

Note that in contrast to the boundary element method the global matrices in FE and IGA are usually 
very sparse due to the compact support of the shape functions. An entry in $M_{ab}$ is only non-zero if the associated shape functions  
$N_a$   and $N_b$ share at least one element where both functions are non-zero.

\subsection{Finite elements}
In the case of finite elements the  shape functions are often simple piecewise polynomials that are continuous at element interfaces.
\begin{figure}[htbp]
\begin{center}
\includegraphics[width=0.45\textwidth]{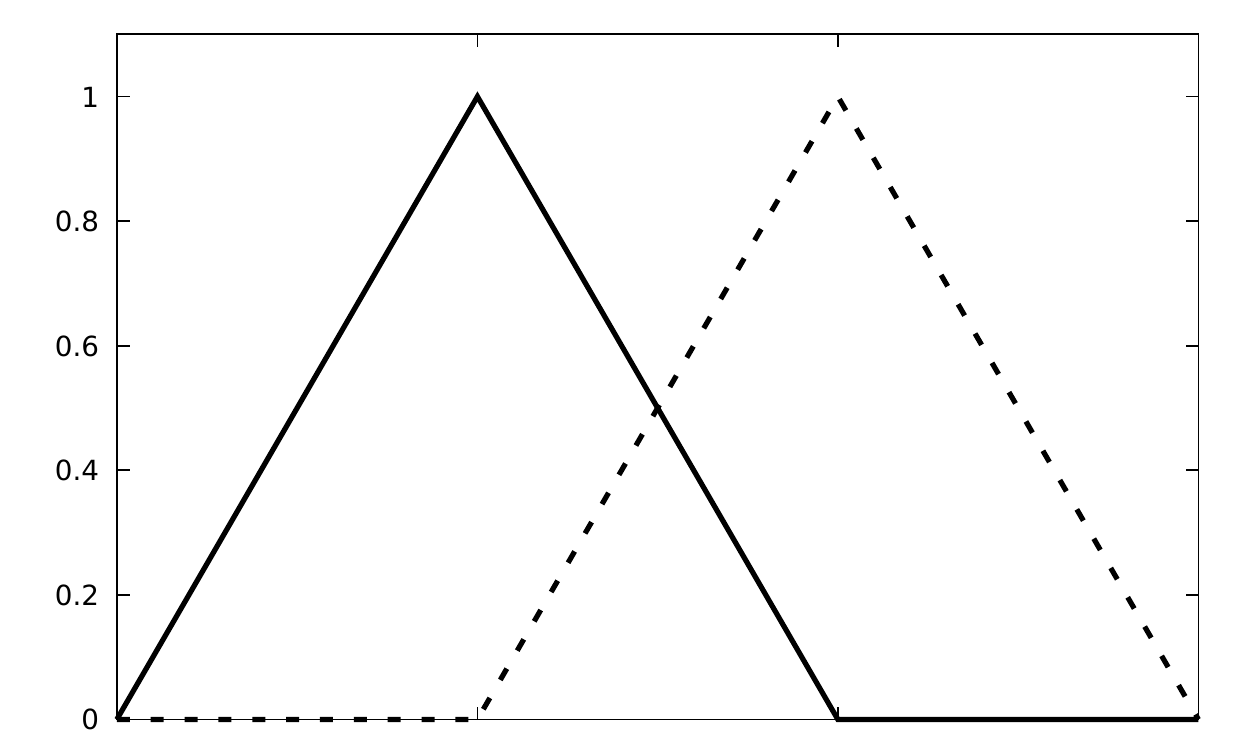}
\caption{Linear shape functions}
\label{fig:lin}
\end{center}
\end{figure}
The simplest example being the 1D linear shape functions as depicted in figure \ref{fig:lin}.

\begin{figure}[htbp]
\begin{center}
\includegraphics[width=0.45\textwidth]{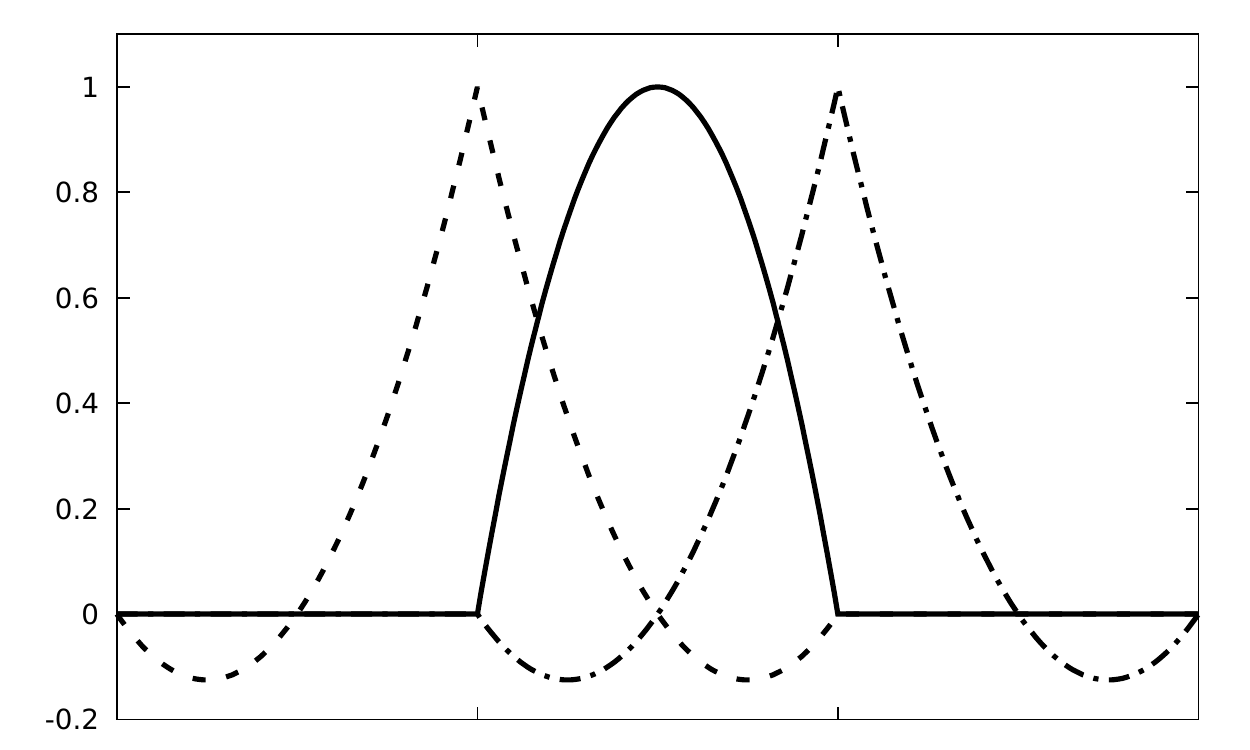}
\caption{Quadratic Lagrangian shape functions}
\label{fig:quad}
\end{center}
\end{figure}

Higher-order shape functions, such as the quadratic Lagrangian shape functions depicted in figure \ref{fig:quad}, are also available.
The expected convergence of the overall numerical method is often directly related to the polynomial order of the shapefunctions used.
Note that the functions are still continuous at element interfaces, but that higher order derivatives are not.

%=============================================
\subsection{Isogeometric analysis}\label{sec:iga}
\begin{figure}[!ht]
\begin{center}
\includegraphics[width=0.45\textwidth]{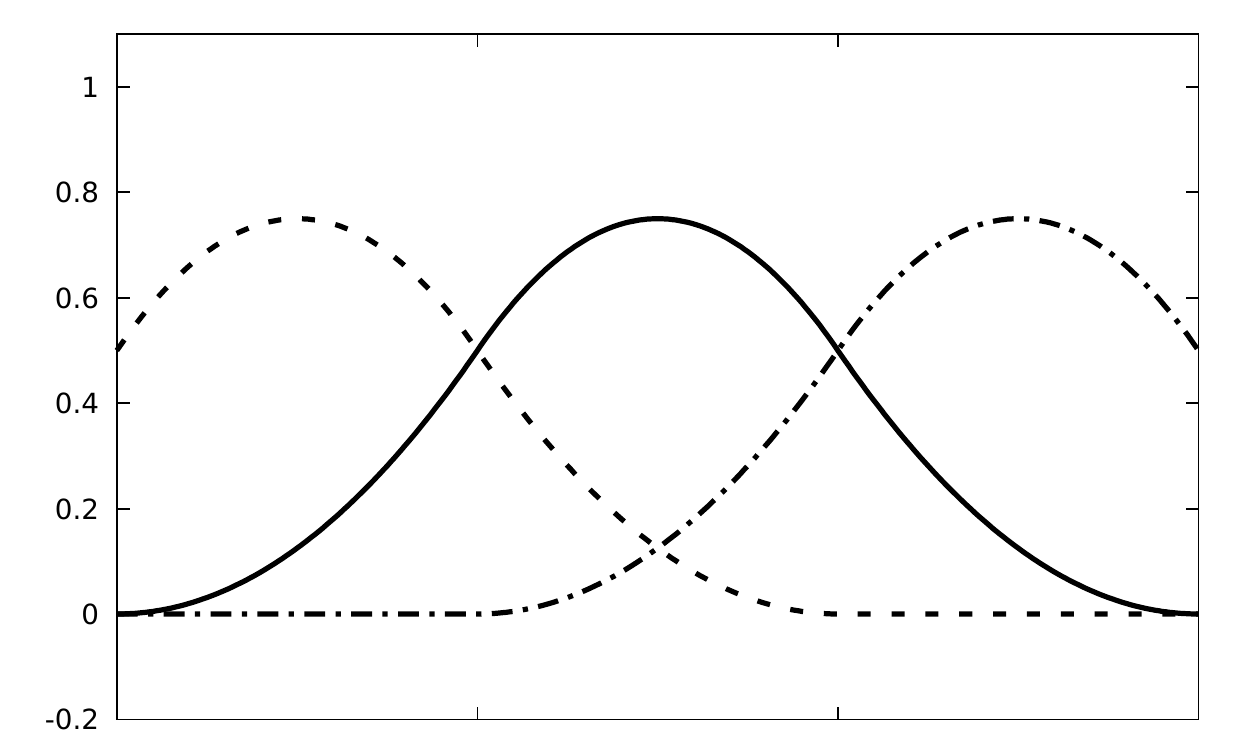}
\caption{Quadratic NURBS shape functions}
\label{fig:nurbs}
\end{center}
\end{figure}

%=============================================
Isogeometric analysis is an extension of finite elements attempting to bridge the gap between design (CAD) and analysis (FE). 
The idea is to use the same NURBS (Non-Uniform Rational B-Splines) of CAD for the analysis.
This results in an exact representation of CAD geometries, such as circles and ellipses. 
These shape functions allow for higher-order continuity at element interfaces, see Figure \ref{fig:nurbs}.

The higher-order continuity leads to additional constraints, resulting in higher-order approximations with limited degrees-of-freedom.
Additionally, the higher-order continuity allows a strict tailoring of approximations spaces.
For the incompressible Navier-Stokes equations this can be exploited to create velocity-pressure approximations that result in exactly solenoidal solutions \cite{Evans13unsteadyNS}.
In \cite{EiAk17ii,AkEik19} this exact divergence is essential for getting correct energy behavior of the single and two-fluid Navier-Stokes problem, respectively.  
In the current context the improved approximation and spectral properties  are of paramount importance. 

In \cite{BACHH07} it is demonstrated that the extra resolution provided by the element kinks are not beneficial.
Evans et al. \cite{EBBH09} show that these kinks lead to bad performance of higher-order finite elements in wave propagation problems.
Cotrell et al. \cite{CoHuRe07} formally prove that NURBS have optimal approximation properties.
This superior behavior is noticed in numerous application areas.

Due to the success of IGA on the one hand and the specific requirements on the other hand, a multitude of alternative spline technologies have emerged. 
Examples are T-Splines \cite{Tspline2010}, LR-Splines \cite{LRspline2013}, U-Splines \cite{USpline2019} and many others.
These approaches all allow higher-order continuity at element interfaces but avoid the rigid tensorial construction required for NURBS.
Assuring the shape functions are linearly independent and the resulting system matrices are solvable is one of the main issues of these alternative splines.

%=============================================
\section{Numerical comparison of the formulations}
\label{sec:numres}
%=============================================
In this section the different weak formulations are numerically investigated.
To this purpose we present a traveling wave case.
We compare the energy behavior and assess convergence under mesh refinement.
%Using the monolithic formulation we assess the relative performance of the isogeometric analysis versus finite elements.
%A validation study is performed by reproducing the dispersion relation.

\subsection{Two-dimensional traveling wave}
The performance of the different formulations is assesses using a simple traveling wave in 
a two-dimensional periodic domain. 

\begin{figure}[!ht]
\begin{center}
\includegraphics[width=0.4\textwidth]{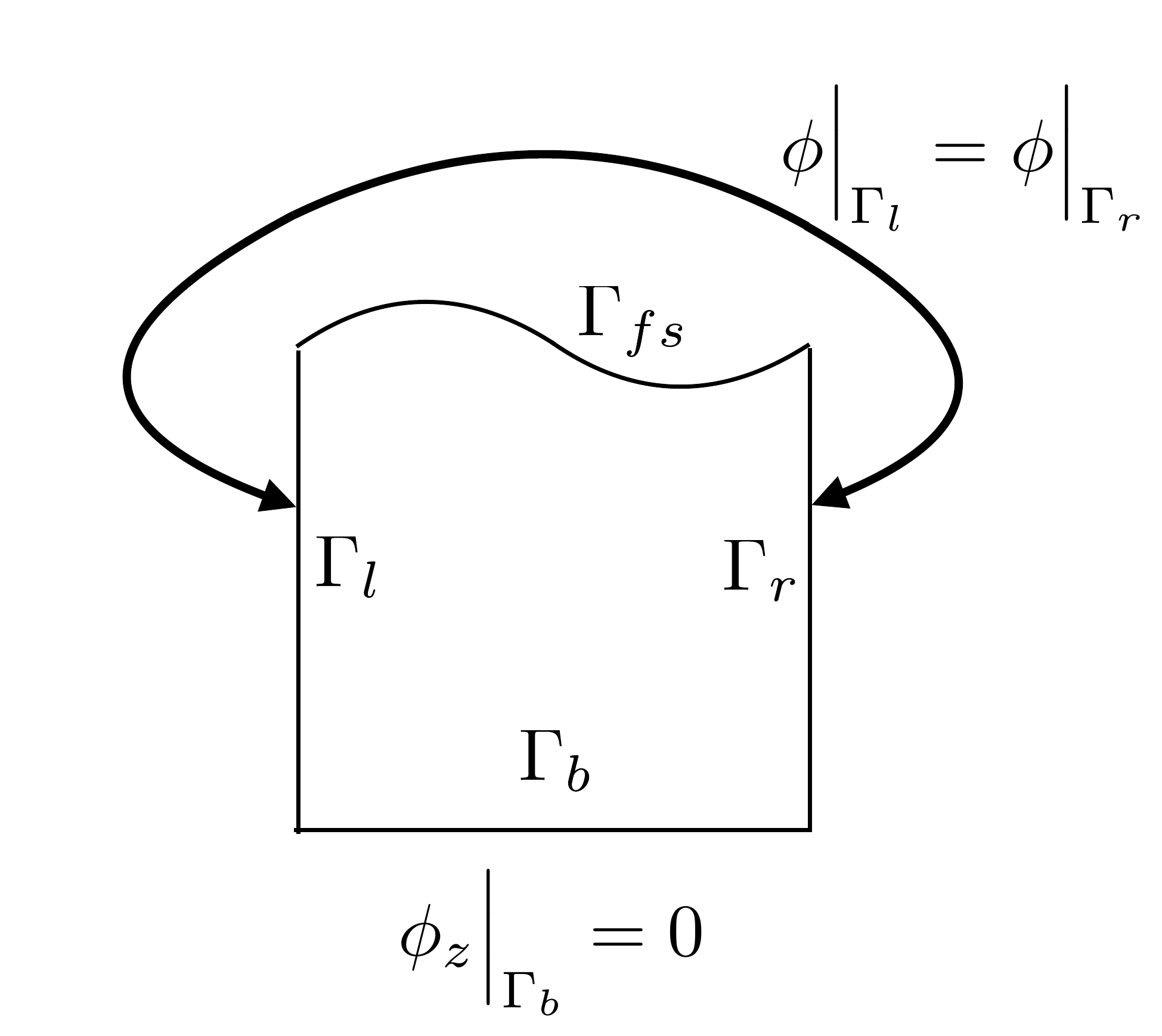}
\caption{Problem setup for simple traveling wave}
\label{fig:setup2D}
\end{center}
\end{figure}
%The domain is a square with depth and length of $1.0m$. 
 Periodic boundary conditions 
 are enforced on the sides ($\Gamma_l$ and $\Gamma_r$), no penetration on the bottom ($\Gamma_b$)
and a free-surface on the top ($\Gammafs$), as depicted in figure \ref{fig:setup2D}.

\begin{figure}[!ht]
\begin{center}
\includegraphics[width=0.15\textwidth]{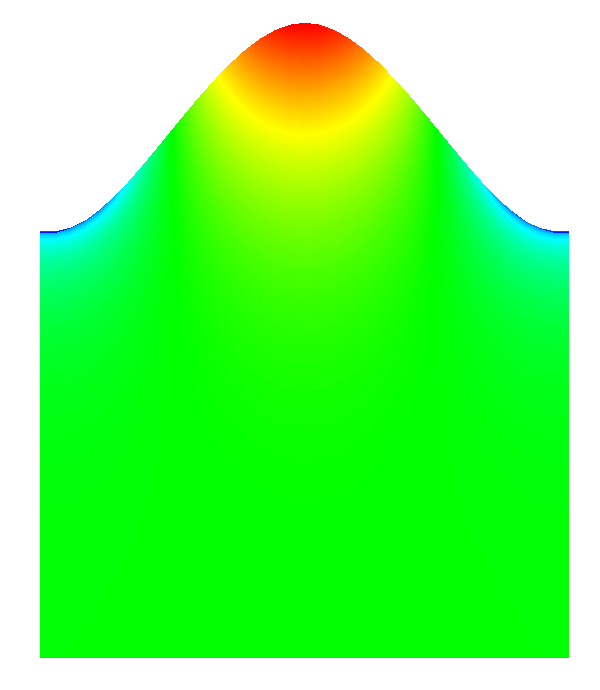}
\includegraphics[width=0.15\textwidth]{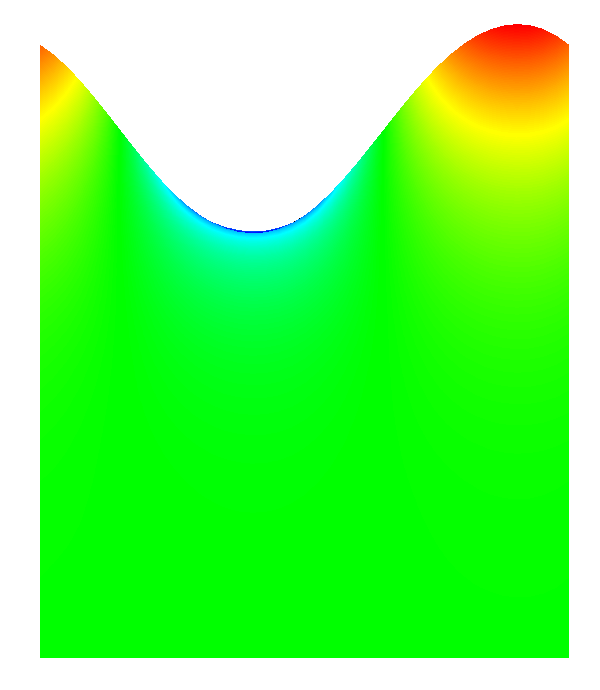}
\includegraphics[width=0.15\textwidth]{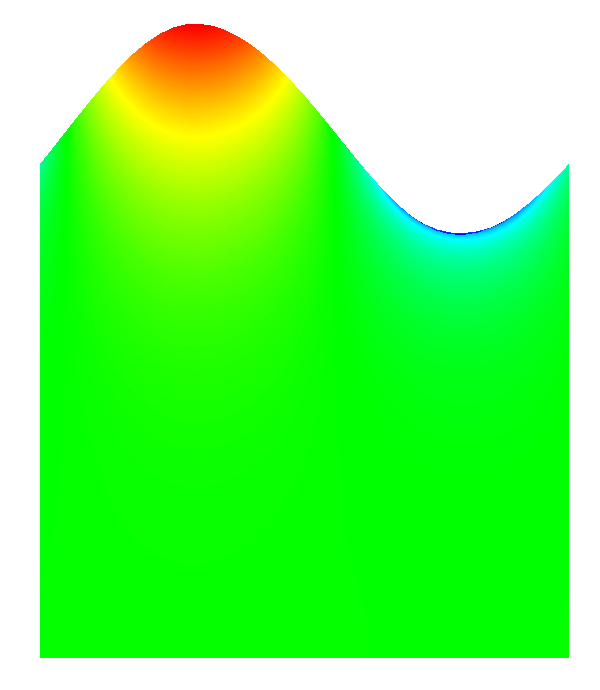}
\caption{Snapshots of the traveling wave problem}
\label{fig:snap2D}
\end{center}
\end{figure}

The initial condition is specified to be an airy wave. For a given wave height $\xi$ the water elevation and flow potential are given as:
\begin{subequations}
\begin{alignat}{1}
\eta \,=~&\xi\cos \,(kx\,-\,\omega t), \\
\phi \,=~&\,{\frac  {\omega }{k}}\,\xi\,{\frac  {\cosh \,{\bigl (}k\,(z+H){\bigr )}}{\sinh \,(k\,H)}}\,\sin \,(kx\,-\,\omega t),
\end{alignat}
\end{subequations}
with $x$ and $z$ representing the spatial coordinates in horizontal and vertical direction.
The wavenumber $k = 2\pi/\lambda$ with wavelength $\lambda$, and angular frequency $\omega$ are related through the dispersion relation:
\begin{align}
\omega ^{2}=g\,k\,\tanh(k\,H),
\end{align}
where $H$ is the water depth.
Snapshots for the solution of an airy-wave with a wavelength equal to the size of the domain are given in Figure \ref{fig:snap2D}.

\clearpage
\input{\figopt/trace2_p1.tex}

\input{\figopt/trace2_p2.tex}
\clearpage
\subsection{Mesh convergence of energy traces}

The simulation are done on a  unit size square. The time step is chosen in alignment with the mesh resolution:
\begin{align}
\frac{h}{\lambda} = \frac{10 \Delta t}{T_{ex}},
\end{align}
where $T_{ex}$ is the theoretically expected period of the wave.

A straightforward computation reveals that the relations for the wave energy in this domain are: 
\begin{subequations}
\begin{alignat}{1} \label{eq:kinpot}
E_{\text{kin}} &=  \frac{1}{2} \|\nabla \phi\|^2_\Omega= \onequart g \xi^2,\\
E_{\text{pot}} &=  \frac{g}{2} \|\eta\|^2_\Gammafs = \onequart g \xi^2, \\
E_{\text{tot}} &=  E_{\text{kin}} +  E_{\text{pot}} = \onehalf g \xi^2. 
\end{alignat}
\end{subequations}
%The normalized time traces of these energies are given in the figures 
%\ref{fig:tracef3p1},
%\ref{fig:tracef2p1},
%\ref{fig:tracef3p2} and
%\ref{fig:tracef2p2}.

Figures \ref{fig:tracef3p1} and \ref{fig:tracef2p1} show the convergence of the energy time trace for linear finite elements 
when employing the segregated and monolithic formulation. 

The energy traces for the reduced formulation are virtually identical to those of the monolithic formulation and are therefore not plotted.
The segregated formulation displays quite severe fluctuations in the energy components and the sum of both does not remain constant. 
These fluctuations disappear with mesh refinement indicating the method is in principle valid.
The monolithic formulation gives much better results on the same mesh when compared with the segregated formulation.  
The fluctuations of the kinetic and potential energy are significantly less than for the segregated formulation.
On the 12 x 12 mesh the fluctuations are barely visible and on the 24 x 24  mesh they are essentially gone.
Moreover, the total energy stays perfectly conserved, even on the coarsest mesh.
Note that on the coarser meshes the total energy is underestimated.
This is due to the loss of energy in projecting the initial condition on to the discrete space.
This mismatch disappears under mesh refinement.

Figures \ref{fig:tracef3p2} and \ref{fig:tracef2p2} show a similar comparison but for quadratic finite elements.
To account for the additional degrees-of-freedom the timestep half compared with the linear cases.
The results are much better owing to the improved resolution from the quadratic elements, even for similar degree-of-freedom counts.
The fluctuations of kinetic and potential energy are significantly reduced, for the monolithic formulation they are 
virtually non-existent even on the coarsest mesh.  
As a consequence, the total energy mismatch for the monolithic formulation is also virtually absent.
The monolithic formulation already has a converged energy behavior on the coarsest mesh of only 3 x 3 quadratic elements.

\subsection{Verification: Mesh convergence of period}
Here we focus on the prediction of the period of the wave.
In practice this can be of major importance as it determines when wave groups or wave crests arrive at certain points.
In the study of wave interference the exact timing of the arrival of each wave becomes critical.

Here we use the same setup as in the previous section.
The simulation is performed for 10 periods. 
The first 2 zero crossings are discarded and the period is computed based on the following 18 zero crossings.

\begin{figure}[htbp]
\begin{center}
\begin{subfigure}{0.47\textwidth}
\centering
\includegraphics[width=\textwidth]{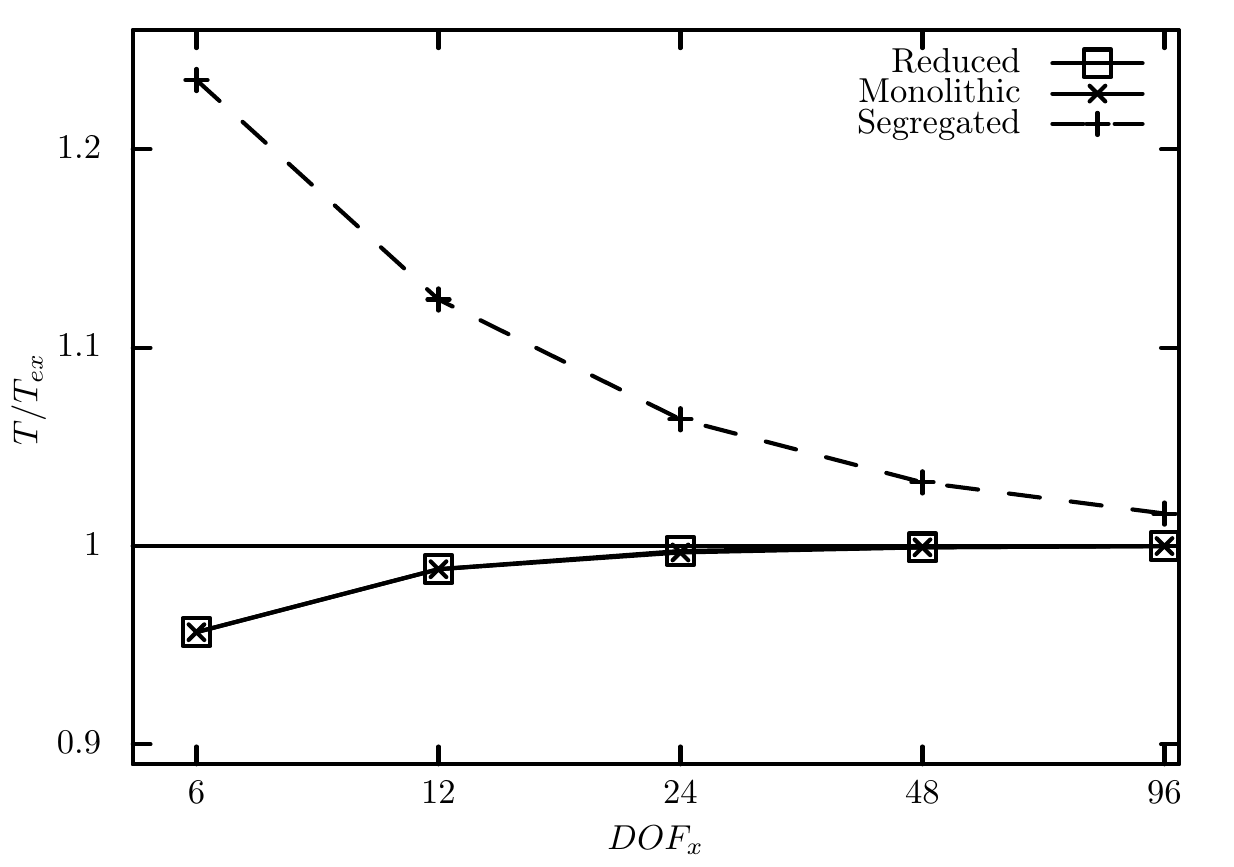}
\caption{Linear finite element}\label{fig:period_lin}
\end{subfigure}
\begin{subfigure}[b]{0.47\textwidth}
\centering
\includegraphics[width=\textwidth]{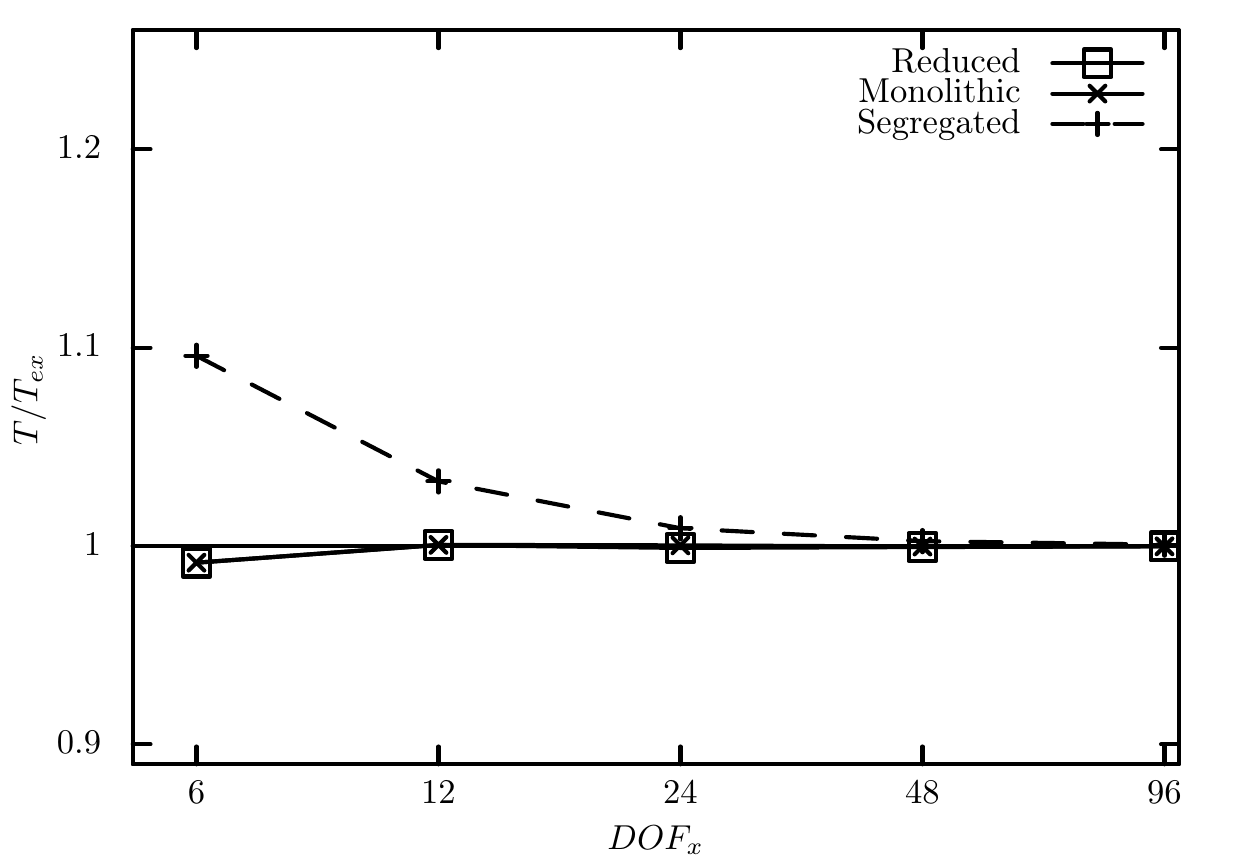}
\caption{Quadratic finite element}\label{fig:period_quad}
\end{subfigure}
\end{center}
\caption{Convergence of the period for different formulations}
\label{fig:period}
\end{figure}

In figure \ref{fig:period} we present the convergence for the different formulations for linear finite elements (figure \ref{fig:period_lin}) and quadratic finite elements (figure \ref{fig:period_quad}).
The reduced and monolithic formulations have virtually the same behavior.
In both cases the segregated formulation performs worse. 
All formulations benefit from the use of quadratic over linear finite elements. 
%However, it is interesting to note that the segregated formulation benefits from the use of quadratic over linear finite elements, 
%while for the other formulations the converse is true. 
%We conjecture the deterioration of the period convergence for  the segregated and monolithic formulations   
%is due to the bad spectral behavior of the higher order finite elements \cite{EBBH09}.
%For the segregated formulation this negative effect is offset by a better prediction of the problematic normal 
%derivative at the free-surface (viz. $\phi_{z}  \vert_\Gammafs $).

\begin{figure}[!ht]
\begin{center}
\includegraphics[width=0.48\textwidth]{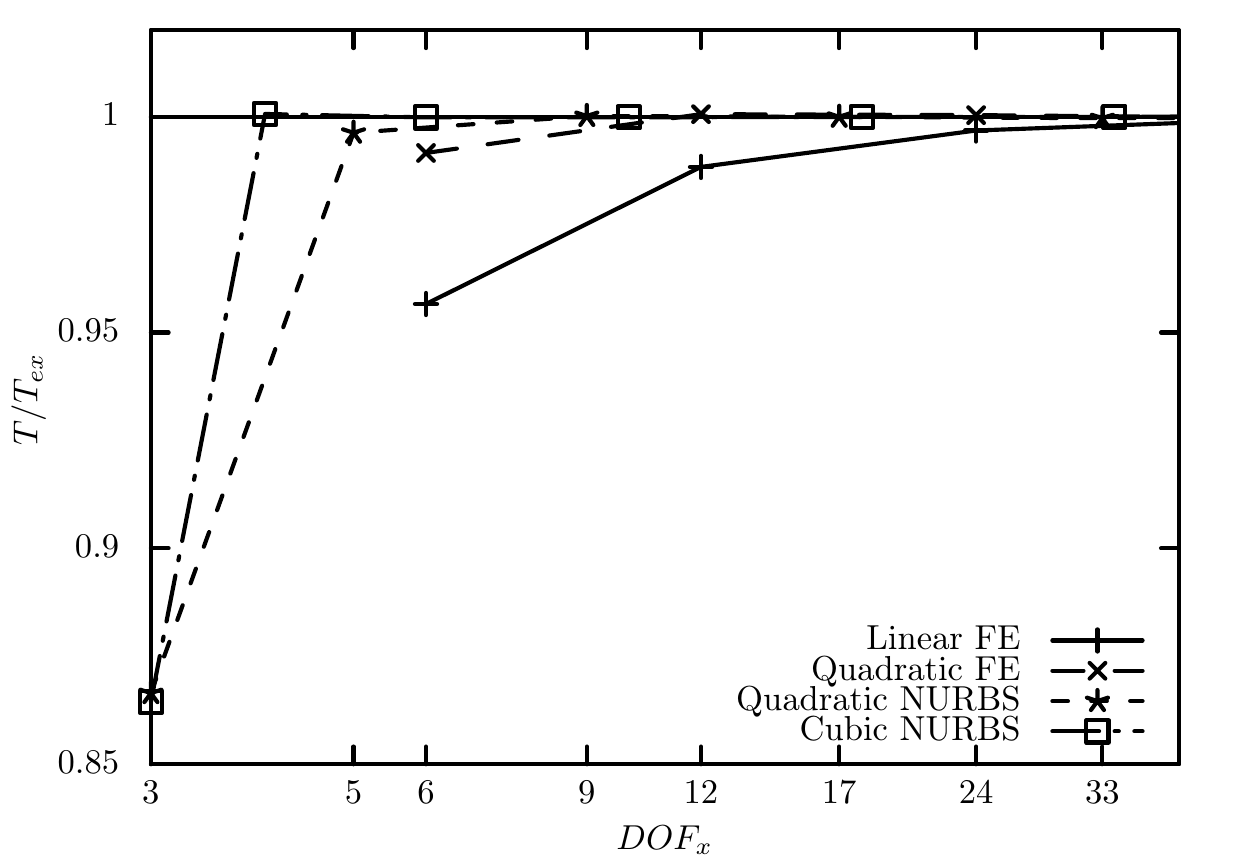}
\end{center}
\caption{Convergence of the period for different discretisations}
\label{fig:period_f2}
\end{figure}

Figure \ref{fig:period_f2} shows the period convergence of the monolithic formulation using different basis functions.
Increasing the order and continuity of the shape functions clearly demonstrates improved behavior.

\begin{figure}[!ht]
\begin{subfigure}[b]{0.23\textwidth}
\begin{center}
\includegraphics[width=\textwidth]{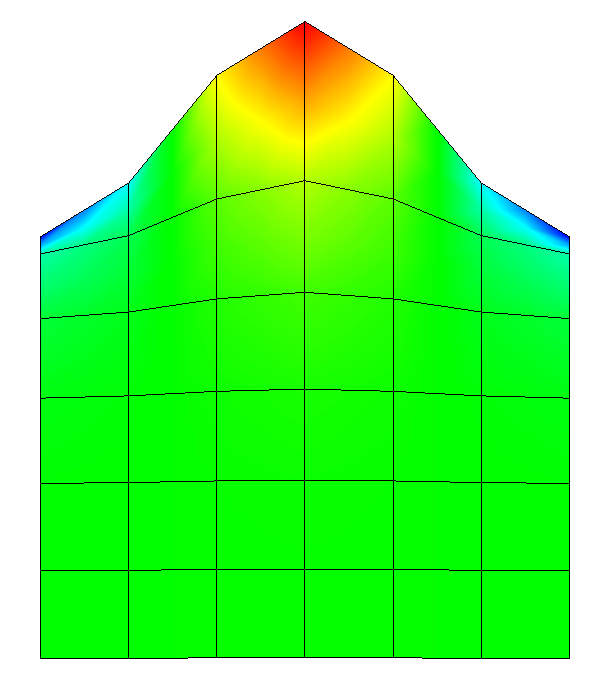}
\caption{Linear FE }
\end{center}
\end{subfigure}
\begin{subfigure}[b]{0.23\textwidth}
\begin{center}
\includegraphics[width=\textwidth]{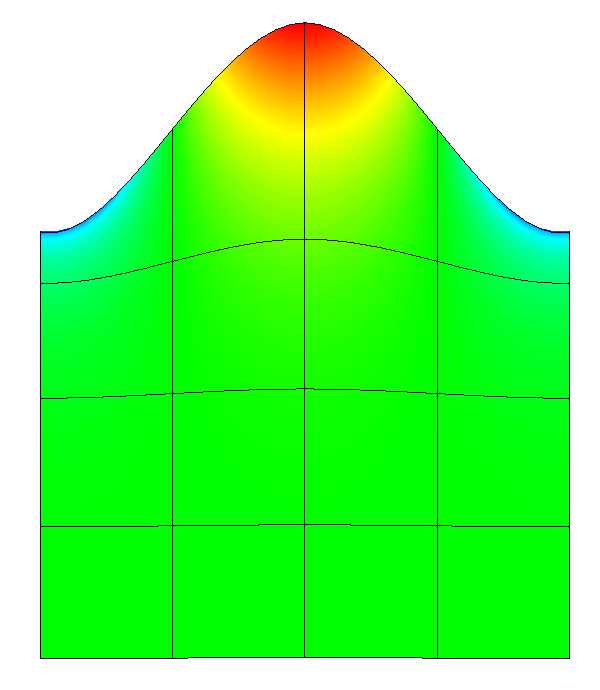}
\caption{Cubic NURBS}
\end{center}
\end{subfigure}
\caption{A visualization of the wave profile for different discretizations.}
\label{fig:ic6}
\end{figure}

Figure \ref{fig:ic6} illustrates the benefit of the higher-order higher-order continuity NURBS basis functions over the standard linear finite elements.
In both cases the horizontal discretization consists of only 6 degrees-of-freedom. 
For the linear  finite element case this leads to a rough estimate and it results in a bad approximation of the period.
In case of cubic NURBS, the spatial representation seems nearly perfect.
The consequence is that the period is predicted very accurately (an error of less than $0.015\%$). 

Given that the monolithic and reduced formulations are nearly indistinguishable and show superior performance, 
we will focus on the monolithic formulation in the following.

%=============================================
\section{Analysis of the weak formulations }
\label{sec:coerc}
%=============================================
In this section we analyze the existence, uniqueness and accuracy of the solution for the monolithic weak formulation.
For completeness the analysis for the reduced and staggered formulations are given in \ref{app:coerc_red}
and \ref{app:coerc_seg}, respectively.
Before we perform the analysis we will shortly introduce the required mathematical tools.

%=============================================
\subsection {Preliminaries}
\label{sec:prelim}
%=============================================
To facilitate the analysis  of the existence, uniqueness and accuracy of the approximate solutions 
we first introduce some essential technical concepts in this section.
Existence and uniqueness of the solutions of the weak formulations are guaranteed  by the Lax-Milgram theorem.
See for instance \cite{ErnGuermondFEM2004} or other standard works on finite element theory. 
\begin{thm}[Lax-Milgram]\label{thm:LM}
Let $V$ be a Hilbert space, $V'$ its dual space, 
and ${B}(\cdot,\cdot)$
a bilinear form on $V$ that is both bounded and coercive:
\begin{subequations}\label{eq:coercive and bnd}
\begin{alignat}{1}
{B}(v, \psi)\leq&~ C_{b}\II \psi \II \; \II v\II, \\
{B}(\psi,\psi)\geq&~ C_{c}\II \psi \II^{2},
\end{alignat}
\end{subequations}
where $\II \cdot \II$ is norm on $V$ and $C_{b}, C_c$ are positive scalars.
Then, for any $f\in V'$, there is a unique solution 
$\psi \in V$ to the equation
\begin{align}
{B}(v,\psi)=f(v)\qquad \text{for all} ~~ v\in V.
\end{align}
The solution $\psi \in V$ satisfies the \textit{a-priori} estimate:
\begin{align}
\II \psi \II \leq {\frac  {1}{C_c}} \sup_{v \in V} \frac{| f(v) |}{\II v \II}.
%\IIIu\III \leq {\frac  {1}{c}} |\!| f |\!|_{{V'}}. \nonumber
\end{align}
\end{thm}

The weak formulations can be converted into a semi-discrete formulation by straightforwardly  
approximating the infinite-dimensional function space by a conforming subspace. 
This is called the Galerkin method.
Existence and uniqueness of solutions obtained with the Galerkin method are automatically inherited from the continuous formulations.
The approximate error of the Galerkin solution is given by  Cea's lemma. 

\begin{lem}[Cea]\label{lem: Cea} Given a bounded and coercive bilinear operator $B(\cdot,\cdot)$
and linear operator $f\in V'$, and finite dimensional space $V^h$ approximating $V$.
There is a unique solution $\psi^h\in V^h$ to the equation
\begin{align}
{B}(v^h,\psi^h)=f(v^h)\qquad \text{for all} ~~ v^h\in V^h.
\end{align}
The error of the approximate solution is bounded:
\begin{align}
\II \psi-\psi^h \II \leq \frac{C_b}{C_c} \inf_{v^h \in V^h}\II \psi-v^h\II  \qquad \text{for all} ~ \psi^h \in V^h,
\end{align}
where $C_b$ and $C_c$ are the constants in (\ref{eq:coercive and bnd}).
\end{lem}

Let us assume that we have a decomposition of the domain into a mesh ${M}=\left\{\Omega_e\right\}_{e=1}^{n_{el}}$ with $\Omega = \cup_e \Omega_e$.
%The free-surface boundary $\Gammafs$ is partitioned as $\Gammafs = \cup_{f_{fs}} \Gamma_{f_{fs}}$.
The element size is denote by $h_e$ with maximum $h=\max_e h_e$.
The boundary $\Gamma$ is partitioned as $\Gamma = \cup_{f} \Gamma_{f}$.
The boundary element size is denote by $h_f$ with maximum $h_b=\max_f h_f$.
In the finite element case the decomposition represents the finite element mesh, whereas in the IGA case it represents the NURBS mesh.

Assuming sufficient regularity of the solution $\psi$, a classical convergence analysis provides the {\em a-priori} interpolation estimates:
%\begin{subequations}
%\begin{alignat}{1}
%\|\nabla \psi - \nabla \tilde{\psi}^h\|_{\Omega} \leq&~ C_1(\ psi) h^{p},\\
%\|\psi - \tilde{\psi}^h\|_{\Gamma} \leq&~ C_0(\psi) h^{p-\onehalf},\\
%\|\psi -  \tilde{\psi}^h\|_{\Gamma} \leq&~ C_{\Gamma}(\psi) h^{p},
%\end{alignat}
%\end{subequations}
\begin{subequations} \label{eq:interpolate}
\begin{alignat}{1}
\inf_{v^h \in V^h} \|\nabla \psi - \nabla v^h\|_{\Omega} \leq&~ C_{\Omega}h^{p} \|  \psi \|_{p,\Omega},\\
\inf_{v^h \in V^h}\|\psi - v^h\|_{\Gamma} \leq&~ C_{T} h^{p+\onehalf} \| \psi \|_{p,\Omega},\\
\inf_{v^h \in V^h}\|\psi -  v^h\|_{\Gamma} \leq&~ C_{\Gamma}h^{p+1} \| \psi \|_{p,\Gamma},
\end{alignat}
\end{subequations}
where $C_{\Omega}$, $C_{T}$  and $C_{\Gamma}$ are constants,  $p$ is the minimum degree of the shape functions and 
$ \|  \cdot \|_{p,\Omega}$ and $\| \cdot \|_{p,\Gamma}$ are the norms of the $p^{th}$ derivative over $\Omega$ and $\Gamma$, respectively.
Note that for simplicity we assume $h \geq h_b$, the validity of this assumption depends on the definition of the element and face sizes $h_e$ and $h_f$. 
For more details see for instance \cite{Ciarlet78,BBVCHS06}.

%=============================================
\subsection{Time-discrete monolithic weak formulation}
%=============================================
Let the test function pair be $\bW:=(w,v) \in \VV$ and the trial function pair denote $\bvphi:=(\phi,\eta) \in \VV$.
The  time-discrete weak formulation then becomes:\\

\textit{Find $\bvphi^{n+1/2} \in \VV$ such that for all $ \bW \in \VV$:}
\begin{align}\label{eq: mono formulation}
&(\nabla w,\nabla \phi^{n+1/2})  - (w, \eta_t^{n+1/2})_{\Gammafs}  \nonumber \\ 
&+\onehalf (v + \frac{\alpha}{g} w,\phi_{t}^{n+1/2}  + g \eta^{n+1/2})_{\Gammafs} =0.
\end{align}
%To obtain the time-discrete version we approximate $\phi_{tt} $ by $\phi_{tt} ^{n+1/2}$. 
We combine the relations (\ref{eq:2nd-int})-(\ref{eq:2nd-kin}) to arrive at:
\begin{subequations}\label{eq:1st-onehalf}
\begin{alignat}{1}
\phi_{t}^{n+1/2} =& \frac{2}{\Delta t} \left (\phi^{n+1/2} -\phi^{n}  \right ), \\
\eta_{t}^{n+1/2} =& \frac{2}{\Delta t} \left (\eta^{n+1/2} -\eta^{n}  \right ).
\end{alignat}
\end{subequations}
Employing these relations we arrive at the \textit{time-discrete} problem:\\
\textit{Given $\phi^n, \eta^n$, find $\bvphi^{n+1/2} \in \VV$ such that for all $ \bW \in \VV$:}
\begin{subequations}
\begin{alignat}{1}
 {B}_m(\bW , \bvphi^{n+1/2} ) =&~ {F}_m(\bW), \label{eq:td form mono}
\end{alignat}
\textit{where}
\begin{alignat}{1}
 {B}_m(\bW , \bvphi^{n+1/2} ):=&~(\nabla w,\nabla \phi)-\frac{2}{\Delta t} (w,\eta)_{\Gammafs}\\
                                               &\qquad~+\frac{1}{2} \left (v+\frac{\alpha}{g}w, \frac{2}{\Delta t} \phi + g \eta\right )_{\Gammafs},\nn\\                                                
{F}_m(\bW):=&~ -\frac{2}{\Delta t}  \left(w, \eta_n\right)_{\Gammafs}+ \frac{1}{2} \left(v+\frac{\alpha}{g}w, \frac{2}{\Delta t}  \phi_n\right)_{\Gammafs}.
\end{alignat}
\end{subequations}

\subsection{Existence and accuracy of the monolithic form} \label{sec:coerc_mono}
The coercivity estimate is 
\begin{align}
{B}_m(\bW,\bW) %=&(\nabla w,\nabla w)  - (w, c v)_{\Gammafs}  \nonumber \\ &
%+\onehalf c (v, w)_{\Gammafs}    + \onehalf g \|v\|^2_{\Gammafs} \nonumber \\ &
%+\frac{c^2}{2g} \|w\|^2_{\Gammafs}+\onehalf c (w , v )_{\Gammafs}\nonumber \\ 
%=&~\|\nabla w\|^2 - \frac{2}{\Delta t}(w,v)_\Gammafs\nn\\
%&+ \frac{1}{2}\left(v+\frac{\alpha}{g}w, \frac{2}{\Delta t} w+gv\right)_\Gammafs \nn \\
=&~\|\nabla w\|^2  +\frac{ \alpha }{g \Delta t}\|w  \|^2_{\Gammafs}  +\frac{g}{2}\| v\|^2_{\Gammafs}\nn \\
&~  +\left (  \frac{\alpha}{2}
- \frac{1}{\Delta t} \right)(w,v)_\Gammafs
.
\end{align}
For  $\alpha = 2/\Delta t$ the coercivity estimate is sharp:
\begin{align}\label{eq:coerc_mono}
{B}_m(\bW,\bW) = \vertiii{\bW}^2_m, 
\end{align}
where the norm is defined as,
\begin{align}\label{eq:norm_mono}
 \vertiii{\bW}^2_m =&~\|\nabla w\|^2  +\frac{2}{g \Delta t^2}\|w  \|^2_{\Gammafs}  +\frac{g}{2}\| v\|^2_{\Gammafs}.
\end{align}

Next, we consider boundedness. We write:
\begin{align}
{B}_m(\bW; \bvphi) \leq&
\|\nabla w\|\|\nabla \phi\| \nn\\
&  + \frac{g}{2} \|v\|_{\Gammafs}\|  \eta \|_{\Gammafs}  + \frac{2}{\Delta t^2 g} \|w\|_{\Gammafs}\|\phi  \|_{\Gammafs} \nonumber \\ &
+\frac{1}{\Delta t} (\|w\|_{\Gammafs}\| \eta\|_{\Gammafs} +\|\phi\|_{\Gammafs}\| v\|_{\Gammafs}).
\end{align}
By defining the following shorthand notation
\begin{align}\label{eq: mono_constants}
 x_{1} =&~\|\nabla w\|, & y_{1} =&~\|\nabla \phi\|, \nn\\
 x_{2} =&~\frac{\alpha}{\sqrt{2g}}\|w\|_{\Gammafs}, & y_{2} =&~\frac{\alpha}{\sqrt{2g}}\|\phi\|_{\Gammafs}, \nn \\ 
 x_{3} =&~\frac{1}{2}\sqrt{2g}\| v\|_{\Gammafs}, & y_{3} =&~\frac{1}{2}\sqrt{2g}\|\eta\|_{\Gammafs},
\end{align}
we can write,
\begin{align}\label{eq: mono_bnd_x}
{B}_{fs}(\bW ;\bvphi) \leq&~x_{1}y_{1}+x_{2}y_{2}+ x_{3}y_{3}\nn\\
&+x_{2}y_{3} +x_{3}y_{2} =\bx \cdot \B{A} \by.
\end{align}
Here $\B{A}$ is a symmetric matrix:
\begin{align}
\B{A} =\left (
\begin{array}{ccc}
1&0&0\\
0&1&1\\
0&1&1\\
\end{array}\right )
\end{align}
with maximum eigenvalue $\lambda_{\rm max}=2$.
Using this $\lambda_{\rm max}$ and applying Cauchy-Schwarz on (\ref{eq: mono_bnd_x}) we arrive at the following boundedness estimate:
\begin{align}
{B}_{fs}(\bW ;\bvphi) \leq&~ \lambda_{max}\|\bx\|_2  \|\by \|_2 \nn \\%= 2 \|\bx\|_2  \|\by \|_2 
\leq&2 \vertiii{\bW}_m\vertiii{\bvphi}_m,
\end{align}
where we used the identities $\|\bx\|_2 = \vertiii{\bW}_m$ and $\|\by \|_2=\vertiii{\bvphi}_m$.

Using Cea's lemma we arrive at the following, accuracy estimate:
\begin{align}
 \vertiii{\bvphi - \bvphi^h}_m & \leq 
  2 \inf_{\bW^h \in \VV^h} \vertiii{\bvphi-\bW^h}_m.
\end{align} 
Using the definition of the norm (\ref{eq:coerc_mono}) and the interpolation estimates (\ref{eq:interpolate}) this yields:
\begin{comment}
\begin{subequations}
\begin{alignat}{1}
 \vertiii{\bvphi - \bvphi^h}^2_m & \leq 
  4 \inf_{\bW^h \in \VV^h} \vertiii{\bvphi-\bW^h}^2_m.
   \nn \\
  &=  4 \inf_{\bW^h \in \VV^h} \|\nabla \phi-\nabla w^h\|^2  +\frac{\alpha^2}{2g}\|\phi-w^h  \|^2_{\Gammafs}  +\frac{g}{2}\|\eta- v^h\|^2_{\Gammafs} \nn \\
   &\leq4( \|\nabla \phi-\nabla w^h\|^2  +\frac{\alpha^2}{2g}\|\phi-w^h  \|^2_{\Gammafs}  +\frac{g}{2}\|\eta- v^h\|^2_{\Gammafs})\nn \\
   &\leq 4C_{\Omega}^2h^{2p} \|  \phi \|^2_{p,\Omega}  +4\frac{\alpha^2}{2g}C^2_{T} h^{2p+1} \| \phi \|^2_{p,\Omega}  +4\frac{g}{2} C^2_{\Gamma}h^{2p+2} \| \eta \|^2_{p,\Gammafs} \nn\\
      &\leq 4 \left ( 1  +\frac{C^2_{T}}{C_{\Omega}^2} \frac{\alpha^2h }{2g}\right )C_{\Omega}^2h^{2p} \|  \phi \|^2_{p,\Omega}  +2 g C^2_{\Gamma}h^{2p+2} \| \eta \|^2_{p,\Gammafs}
\end{alignat}
\end{subequations} 
\end{comment}
\begin{align}
 \vertiii{\bvphi - \bvphi^h}^2_m 
 & \leq 
4 \left ( 1  +\frac{C^2_{T}}{C_{\Omega}^2} \frac{2h }{g \Delta t ^2}\right )C_{\Omega}^2h^{2p} \|  \phi \|^2_{p,\Omega}  
\nn\\& \qquad\qquad\qquad
+2 g C^2_{\Gamma}h^{2p+2} \| \eta \|^2_{p,\Gammafs}.
\end{align}
This indicates that depending on the relative importance of the terms we could pick an additional full order of convergence order.

%=============================================
\section{V\&V of the monolithic formulation}
\label{sec:VV}
%=============================================

\subsection{Verification: Error convergence}
Here we present the error convergence for the monolithic formulation.
We measure the error in the norm (\ref{eq:coerc_mono}).
Note that this norm includes the timestep size $\Delta t$.
Therefore we choose to perform the mesh convergence study with a fixed timestep. 
In this way the definition of the norm does not change when refining the mesh.
The timestep is chosen as $\Delta t = T/1000$, with end time $T$, leading to negligible time stepping errors.
\begin{figure}[!ht]
\begin{center}
\includegraphics[width=0.4\textwidth]{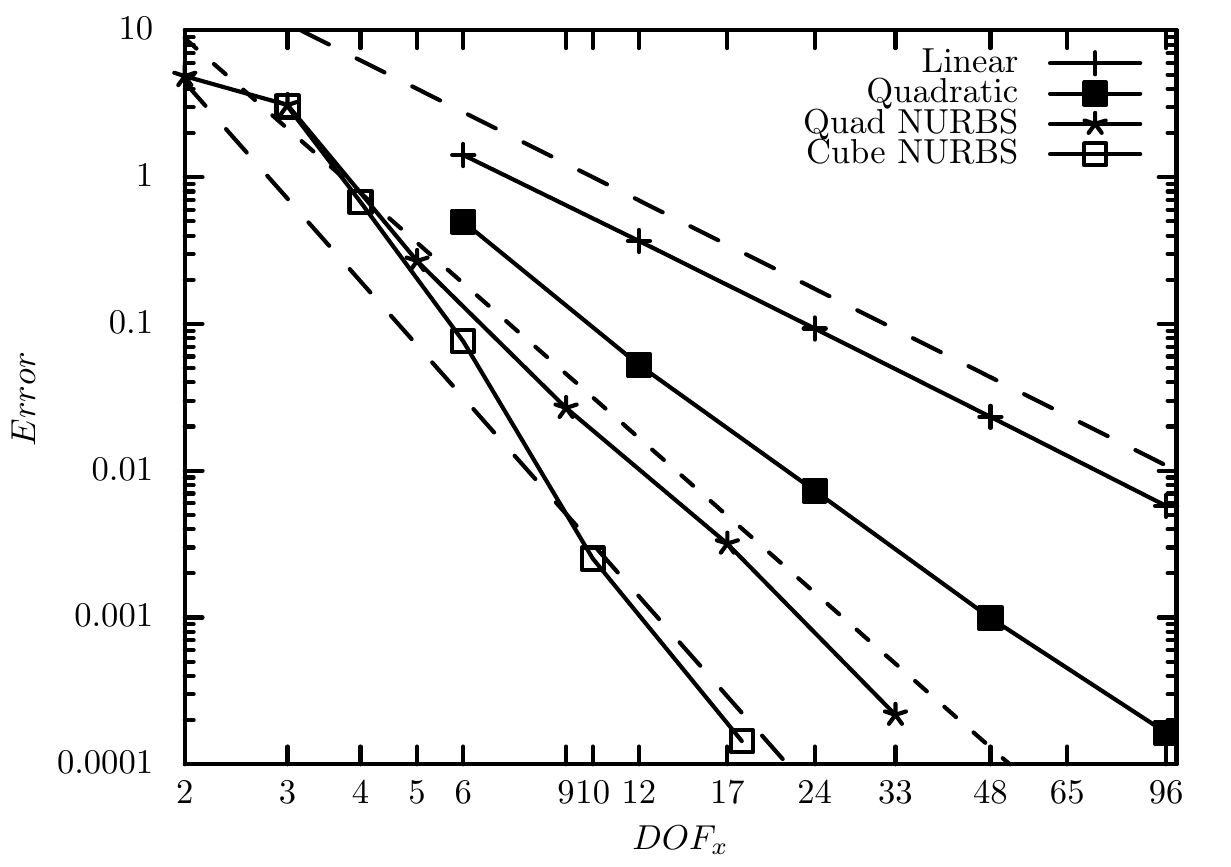}
\end{center}
\caption{Convergence of error for the monolithic formulation. The dashed lines indicate convergence of the orders $O(h^2)$, $O(h^{3\onehalf})$  and $O(h^{4\onehalf})$ respectively.}
\label{fig:error_conv}
\end{figure}
The convergence of the error is shown in figure \ref{fig:error_conv}.
Note that the theoretical convergence rates are $\mathcal{O}(h^p)$ and $\mathcal{O}(h^{p+1})$, with $p$ the order of the basis functions, for the volumetric and boundary term.
This indicates that for all discretizations, except for linear finite elements, the potential boundary term is the dominant term.
Remark that obtained convergence rates are larger than the theoretical values, which confirms the viability of the formulation.

\subsection{Validation: Dispersion relation}

To show the adequacy of the monolithic formulation we will perform a validation case which consists of predicting the phase velocity of the wave for different heights.
The theoretical solution is:
\begin{align}
c_p = \sqrt{\frac{g}{k} \tanh (kH)},
\end{align}
with phase velocity $c_p$.
The numerical approximation is obtained on a mesh with 8 x 8 cubic NURBS elements resulting in  a 10 x 11 degrees-of-freedom resolution.
The time step is taken such that $80$ periods fit in the time domain.

\begin{figure}[!ht]
\begin{center}
\includegraphics[width=0.4\textwidth]{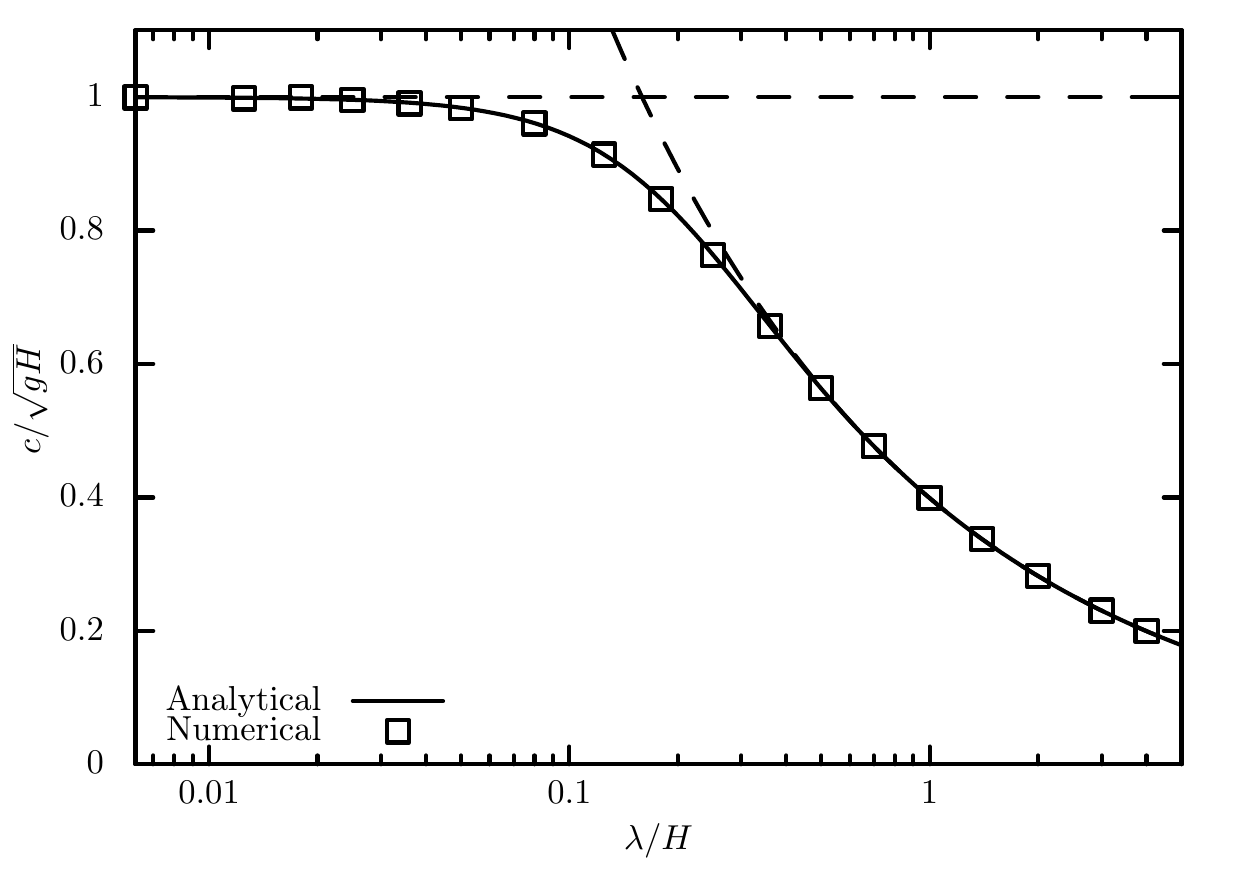}
\caption{Dispersion relation}
\label{fig:VV-dispersion}
\end{center}
\end{figure}

The phase speed is computed from the wave period in the same manner as in the previous subsections.
Figure \ref{fig:VV-dispersion} shows that the numerical method gives a near perfect prediction of the wave 
speed over a wide range of water depths. 
This demonstrates that the monolithic formulation produces the correct physical answer.

%=============================================
\section{Showcase: Sloshing in closed 3D container}
\label{sec:show}
\begin{figure}[!ht]
\begin{center}
\includegraphics[width=0.24\textwidth]{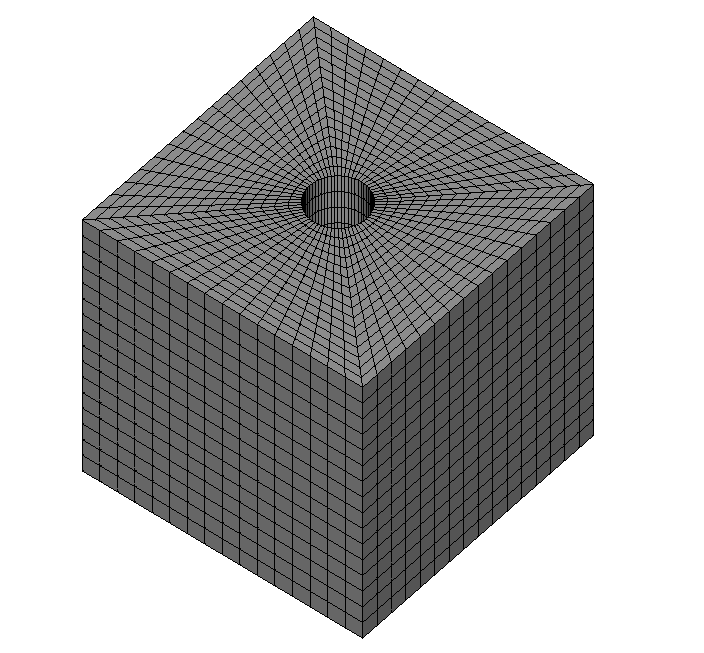}
\caption{Mesh of the 3D sloshing problem}
\label{fig:setup3D}
\end{center}
\end{figure}

\begin{figure}[!ht]
\begin{center}
\includegraphics[width=0.23\textwidth]{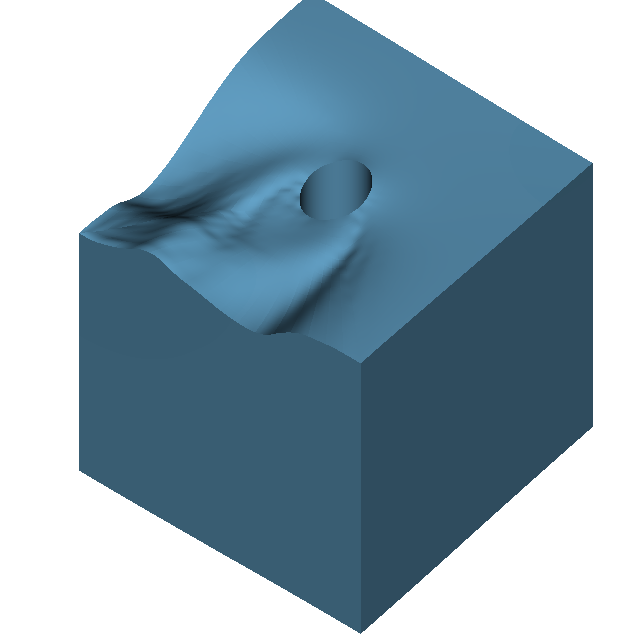}
\includegraphics[width=0.23\textwidth]{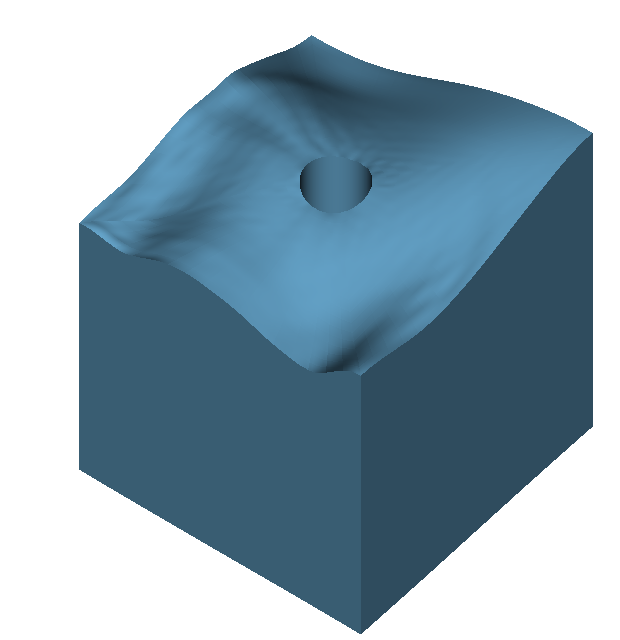}\\
\includegraphics[width=0.23\textwidth]{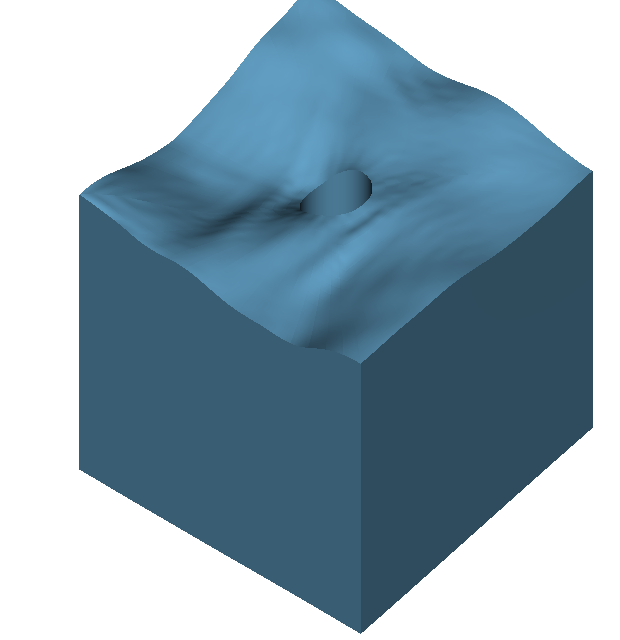}
\includegraphics[width=0.23\textwidth]{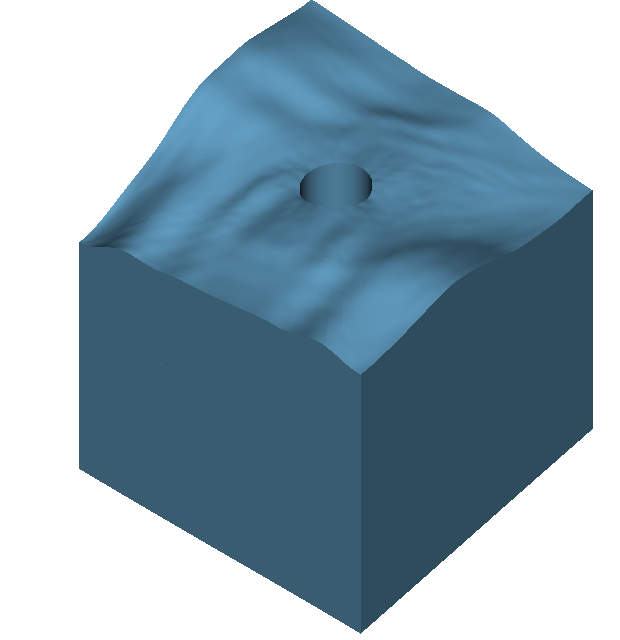}
\caption{Snapshots of standing waves}
\label{fig:snapshots}
\end{center}
\end{figure}

\begin{figure}[!ht]
\begin{center}
\includegraphics[width=0.45\textwidth]{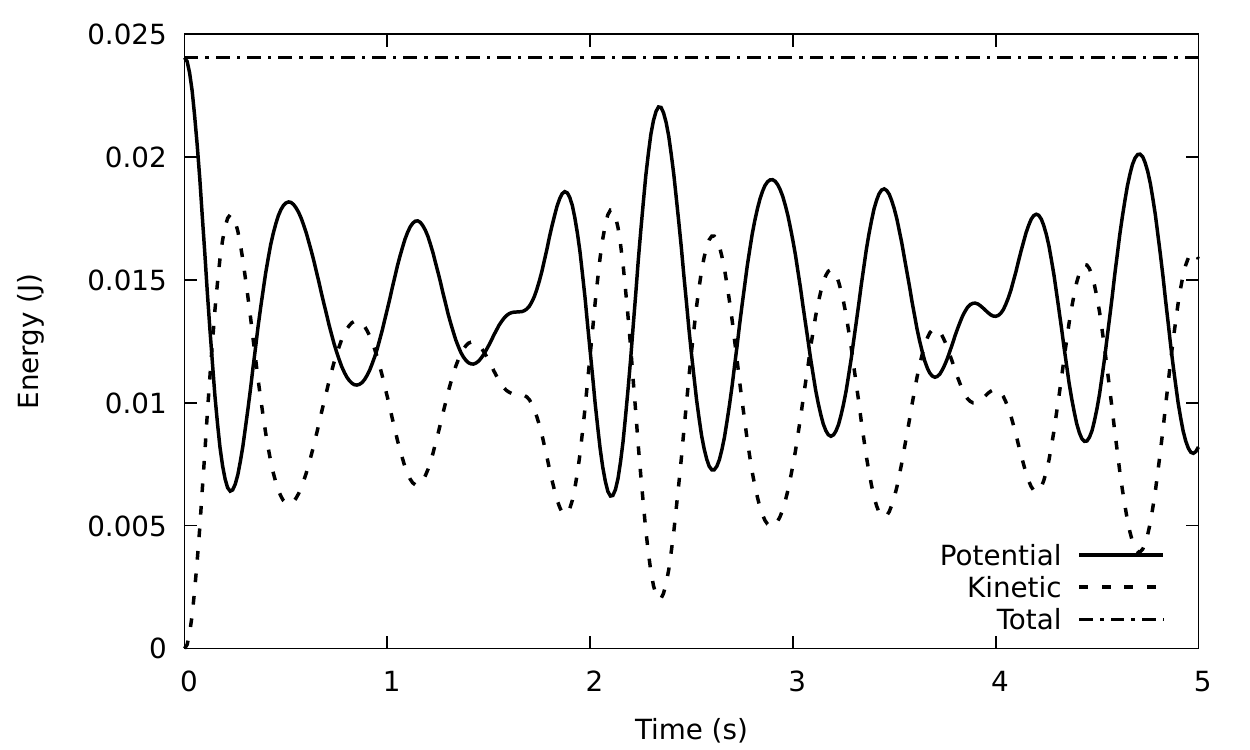}
\caption{Energy trace of the 3D sloshing problem}
\label{fig:trace3D}
\end{center}
\end{figure}

%=============================================
Here we consider a sloshing case in a three-dimensional container with a circular obstruction in the middle of the domain.
The domain is represented with a 32x32x32-mesh with no penetration boundary conditions on all boundaries.
Note that IGA offers the ability to represent circles exactly.
A picture of the setup is presented in figure \ref{fig:setup3D}.
The water surface is given an initial perturbation in one quarter of the domain.
Some snapshots of the resulting water surface using the monolithic formulation are shown in figure \ref{fig:snapshots}.
The resulting traces of kinetic, potential and total energy are depicted in figure \ref{fig:trace3D}.
It shows that kinetic and potential energy are exchanged in an erratic fashion due to the bouncing waves between the cylinder and the outer boundaries.
Just as in the 2D cases the total energy is perfectly conserved.

%=============================================
\section{Conclusion}
\label{sec:conc}
%=============================================
This paper presents a novel variational formulation to simulate linear free-surface flow. 
Similar to formulations found in literature the potential and water height are the primary variables. 
This should facilitate the extension to non-linear cases and situations with forward speed.

The novelty of the formulation is that the interior problem and kinetic and dynamic boundary 
conditions are combined in one formulation. This monolithic formulation results in provable 
energy conservation, if a proper time integrator is selected. 
Additional to the superior energy behavior, the monolithic formulation exhibits significantly reduced dispersion errors.

The dispersion properties of the numerical method are further improved by the adoption of
Isogeometric analysis.
Isogeometric analysis is a discretization approach that at the same time allows for 
higher-order finite elements and higher-order and higher-continuity basis functions. 
This later feature results in more efficient and improved behavior of the numerical solutions.

We demonstrate that the benefit of employing Isogeometric analysis over traditional finite elements is substantial. 
The error in the kinetic and potential energy is very small for quadratic NURBS on very coarse meshes.
Moreover, on very coarse discretizations (4 to 6 degrees-of freedom per wavelength) the dispersion error is negligible.

The analysis of numerical results is performed for a two-dimensional linear wave problem.
A simple sloshing case demonstrates that the results directly translate to three-dimensions.
Possible further work entails the extension to nonlinear cases with forward speed.
%=============================================
\section*{Acknowledgements}

The authors gratefully acknowledge support from Delft University of Technology.
The numerical simulation are performed using the open-source library {MFEM} \cite{mfem-library}, while 
results are visualized using Vis{I}t \cite{VisIt}.% and gnuplot \cite{gnuplot}.

%=============================================
%\section*{References}
\bibliographystyle{unsrt}
\bibliography{references}
%=============================================
\newpage
\appendix

%=============================================
\section{Analysis of the reduced form}\label{app:coerc_red}
%=============================================
To analyze the properties of this variational form, we first discretize in time.
Approximating $\phi$ by $\phi^{n+1/2}$ and $\phi_{tt} $ by $\phi_{tt} ^{n+1/2}$ gives:\\

\textit{Given $\phi^n, \phi_t^n, \phi_{tt} ^n$, find $\phi^{n+1/2} \in \WW$ such that for all $ w \in \WW$:}
\begin{align}\label{eq: reduced weak formulation no substitution}
 (\nabla w,\nabla \phi^{n+1/2})_\Omega  +\frac{1}{g} (w, \phi_{tt} ^{n+1/2})_{\Gammafs} = 0.
\end{align}

In order to prove existence, uniqueness and accuracy estimates, we use Lax-Milgram and Cea's lemma.
To that purpose we establish coercivity and boundedness of the weak form.

We combine the relations (\ref{eq:2nd-int})-(\ref{eq:2nd-kin}) and get:
\begin{align}
\phi_{tt} ^{n+1/2} = \frac{4}{\Delta t^2} \left (\phi^{n+1/2} -\phi^{n}  \right ) - \frac{2}{\Delta}\phi_{t}^{n+1/2}. 
\end{align}
Using this relation we arrive at the following weak form:\\

\textit{Given $\phi_n, \phi_{t}^n, \phi_{tt} ^n$, find $\phi^{n+1/2} \in \WW$ such that for all $ w \in \WW$:}
\begin{subequations}\label{eq:td reduced weak formulation}
\begin{alignat}{1}
 {B}_r(w, \phi^{n+1/2}) =&~ {F}_r(w), \label{eq:td form1}
\end{alignat}
\textit{where}
\begin{alignat}{1}
{B}_r(w,\phi):=&~(\nabla w,\nabla \phi)_\Omega  +\frac{4}{ \Delta t ^2g}  (w, \phi)_{\Gammafs},\\
{F}_r(w):=&~\frac{4}{\Delta t^2 g} (w, \phi_n)_{\Gammafs}+\frac{2}{ \Delta t g} (w, \phi_{t}^n)_{\Gammafs}.
\end{alignat}
\end{subequations}
By defining a problem-dependent norm as 
\begin{align}\label{eq:coerc1}
 \vertiii{w}^2_r := \|\nabla w \|_\Omega^2 +\frac{4}{ \Delta t ^2g} \|w\|^2_{\Gammafs},
\end{align}
the coercivity and boundedness estimates are sharp in this norm:
\begin{subequations}
\begin{alignat}{1}
{B}_r(w,w) =&~ \vertiii{w}_r^2\quad \text{for all }w \in \WW\\
{B}_r(w,\phi) \leq& \|\nabla w\|_\Omega \|\nabla \phi\|_\Omega  +\frac{4}{ \Delta t ^2g}\|w \|_{\Gammafs}  \| \phi\|_{\Gammafs},
  \\  \nonumber
  \leq&  \vertiii{w}_r\vertiii{\phi}_r\quad \text{for all }w \in \WW, \phi \in \WW.
\end{alignat}
\end{subequations} 
The coercivity follows directly from the definition of the norm, while the boundedness estimate 
requires the Cauchy-Schwarz inequality.

Using Cea's lemma we arrive at the following, accuracy estimate:
\begin{align}
 \vertiii{\phi - \phi^h}_r &\leq 
  \inf_{w^h \in W^h} \vertiii{\phi-w^h}_r.
%   \nn \\
%  &=  \inf_{w^h \in W^h} \|\nabla \phi- \nabla w^h \|_\Omega^2 +\frac{4}{ \Delta t ^2g}  \inf_{w^h \in W^h}  \|\phi-w^h\|^2_{\Gammafs}, \nn \\
 %  &\leq \left ( 
% C^2_{\Omega}h
% +   \frac{4C^2_{T}}{ \Delta t ^2g}   \right ) h^{2p-1}  \| \phi \|^2_{p,\Omega}.
\end{align}
Using the definition of the norm (\ref{eq:coerc1}) and the interpolation estimates (\ref{eq:interpolate}) this yields:
\begin{align}
\|\nabla \phi - \nabla\phi^
h \|_\Omega^2 +\frac{4}{ \Delta t ^2g} \|\phi - \phi^h\|^2_{\Gammafs}
 \qquad \qquad  \nn \\ \qquad \qquad 
% \leq \left (  C^2_{\Omega}h+   \frac{4C^2_{T}}{ \Delta t ^2g}   \right ) h^{2p-1}  \| \phi \|^2_{p,\Omega}.
%  \qquad \qquad  \nn \\ \qquad \qquad 
\leq  \left ( 1 +   4 \frac{C^2_{T} }{  C_{\Omega}^2}  \frac{h }{ \Delta t ^2g}   \right ) C_{\Omega}^2 h^{2p}  \| \phi \|^2_{p,\Omega}.
\end{align}
This indicates that depending on the relative importance of the terms we could pick an additional half order of convergence order.

%=============================================
\section{Analysis of the segregated form}\label{app:coerc_seg}
%=============================================

\subsection{Interior problem}

The weak form (\ref{eq:form2a}) is trivially coercive:
\begin{align}
{B}_{int}(w,w)  = \|\nabla w\|_\Omega^2 \quad \text{for all }w \in \WW_0.
\end{align}
Boundedness of the weak form follows directly  from Cauchy-Schwartz:
\begin{align}
{B}_{int}(w,\phi)  \leq &~ \|\nabla w\|_\Omega \|\nabla \phi\|_\Omega \nn\\
&\text{for all }w \in \WW_0, \phi \in \WW_{\hat{\phi}},
\end{align}
in the standard $H^1(\Omega)$-norm.

Using Cea's lemma the accuracy is estimate to be:
\begin{align}
\|\nabla \phi - \nabla \phi^h\|_\Omega &\leq 
  \inf_{w^h \in \WW_0^h}\|\nabla \phi - \nabla w^h\|_\Omega.
\end{align}
Using the interpolation estimates (\ref{eq:interpolate}) this yields:
\begin{align}
\|\nabla \phi - \nabla\phi^h \|_\Omega
\leq  C_{\Omega} h^{p}  \| \phi \|_{p,\Omega}.
\end{align}
This is the standard optimal estimate for standard Poisson problems. 

\subsection{Boundary problem}

%\begin{subequations}\label{eq:1st-onehalf}
%\begin{alignat}{1}
%\phi_{t}^{n+1/2} =& \frac{2}{\Delta t} \left (\phi^{n+1/2} -\phi^{n}  \right ) \\
%\eta_{t}^{n+1/2} =& \frac{2}{\Delta t} \left (\eta^{n+1/2} -\eta^{n}  \right )
%\end{alignat}
%\end{subequations}

To analyze the properties of this variational form, we first discretize in time.
Substitution of the relations (\ref{eq:1st-onehalf}) into (\ref{eq:form2b}) gives:\\

\noindent \textit{Given $\phi^n, \phi_t^n, \eta^n, \eta_t^n$ and $\phi_z$, find $\bvphi^{n+1/2} \in \VV_\Gammafs$ such that for all $ \bW \in \VV_\Gammafs$:}
\begin{subequations}\label{eq:dt form seg b}
\begin{alignat}{1}
{B}_{fs}(\bW ;\bvphi^{n+1/2}) = {F}_{fs}(\bW)
\end{alignat}
\textit{where}
\begin{alignat}{1}
{B}_{fs}(\bW ;\bvphi) :=&~ 
\left (w,\frac{2}{\Delta t} \phi  + g \eta \right )_{\Gammafs}+\frac{g^2}{\alpha^2} \frac{2}{\Delta t}  (v , \eta)_{\Gammafs}, \\
{F}_{fs}(\bW) :=& ~\frac{2}{\Delta t}(w, \phi^n)_{\Gammafs} +  \frac{g^2}{\alpha^2} \frac{2}{\Delta t} (v , \eta^n)_{\Gammafs}\nn\\
&+ \frac{g^2}{\alpha^2} (v ,\phi_{z} )_{\Gammafs}.
\end{alignat}
\end{subequations}
Coercivity of the boundary problem follows via:
\begin{align}
{B}_{fs}(\bW;\bW) =& \frac{2}{\Delta t}\|w\|^2_\Gammafs + g(w,v) +\frac{g^2}{\alpha^2}\frac{2}{\Delta t} \|v\|^2_{\Gammafs} \\
%\geq&~ \frac{2}{\Delta t}\|w\|^2_\Gammafs -\frac{g}{2} \left(\frac{\alpha}{g} \|w\|^2_\Gammafs + \frac{g}{\alpha} \|v\|^2_{\Gammafs} \right)\nn\\
%&~+\frac{g^2}{\alpha^2}\frac{2}{\Delta t} \|v\|^2_{\Gammafs}
 %\nn \\
 %\end{align}
 %\begin{align}
%{B}_{fs}(\bW;\bW) 
   \geq&~ \left ( \frac{2}{\Delta t}-\frac{\alpha}{2}  \right ) \|w\|^2_\Gammafs 
   %\nn\\&~
+\left ( \frac{g^2}{\alpha^2}\frac{2}{\Delta t} -  \frac{g^2}{2 \alpha}   \right ) \|v\|^2_{\Gammafs}.  \nn
% \nn\\
%   =&~ \left ( \frac{2}{\Delta t}-\frac{\alpha}{2}  \right ) \|w\|^2_\Gammafs 
%+\frac{g^2}{ \alpha} \left ( \frac{2}{\alpha \Delta t} -    \onehalf \right ) \|v\|^2_{\Gammafs} 
 \end{align}
Selecting $\alpha = 2/\Delta t$ yields:
\begin{align}\label{eq: norm seg fs}
{B}_{fs}(\bW;\bW)=&~ \vertiii{\bW}_{s}^2,
\end{align}  
with the norm defined as:
\begin{align}
\vertiii{\bW}_{s}^2 :=  \frac{1}{\Delta t} \|w\|^2_\Gammafs + \frac{g^2 \Delta t}{4} \|v\|^2_{\Gammafs}.
 \end{align}
Next, we prove boundedness.
Applying Cauchy-Schwartz we get:
\begin{align}\label{eq:bnd}
{B}_{fs}(\bW ;\bvphi) \leq&~
\frac{2}{\Delta t} \|w\|_{\Gammafs}\|\phi \|_{\Gammafs}
+g \|w\|_{\Gammafs} \|\eta\|_{\Gammafs} \nonumber \\&
+\frac{g^2 \Delta t} {2}\ \|v\|_{\Gammafs} \| \eta\|_{\Gammafs}.
\end{align}
By defining the following shorthand notation:
\begin{align}\label{eq: stag_constants}
 x_{1} &~= \frac{1}{\sqrt{\Delta t}}\|w\|_\Gammafs,& y_{1} &~= \frac{1}{\sqrt{\Delta t}}\|\phi\|_\Gammafs, \nn \\
 x_{2} &~= \frac{ g \sqrt{\Delta t}}{2}\|v \|_\Gammafs,& y_{2} &~= \frac{ g \sqrt{\Delta t}}{2}\|\eta\|_\Gammafs,
\end{align}
we can write,
\begin{align}\label{eq: seg_bnd_x}
{B}_{fs}(\bW ;\bvphi) \leq&~2 x_{1} y_{1} +2 x_{2}y_{2} + 2 x_{1} y_{2} =\bx \cdot \B{A} \by.
\end{align}
Here $\B{A}$ is a symmetric matrix:
\begin{align}
\B{A} =\left (
\begin{array}{cc}
2&2\\
0&2\\
\end{array}\right )
\end{align}
with  maximum eigenvalue $\lambda_{\rm max}=2$.
Using this $\lambda_{\rm max}$ and applying Cauchy-Schwarz on (\ref{eq: seg_bnd_x}) we arrive at: % the following boundedness estimate:
\begin{align}
{B}_{fs}(\bW ;\bvphi) \leq&~ \lambda_{\rm max}\|\bx\|_ \|\by \|_2 \nn \\%= 2 \|\bx\|_2  \|\by \|_2 
\leq& 2\vertiii{\bW}_s\vertiii{\bvphi}_s,
\end{align}
where we used the identities $\|\bx\|_2 = \vertiii{\bW}_s$ and $\|\by \|_2=\vertiii{\bvphi}_s$.

Using Cea's lemma the accuracy is estimated as:
\begin{align}
 \vertiii{\bvphi - \bvphi^h}_s & \leq 
  2 \inf_{\bW^h \in \VV^h} \vertiii{\bvphi-\bW^h}_s.
\end{align} 
Using the definition of the norm (\ref{eq: norm seg fs}) and the interpolation estimates (\ref{eq:interpolate}) this yields:
\begin{comment}
\begin{subequations}
\begin{alignat}{1}
 \vertiii{\bvphi - \bvphi^h}^2_s & \leq 
  2 \inf_{\bW^h \in \VV^h} \vertiii{\bvphi-\bW^h}^2_s.
   \nn \\
&  \leq  inf \frac{2}{\Delta t} \|\phi -w^h\|^2_\Gammafs + \frac{g^2 \Delta t}{2} \|\eta-vh\|^2_{\Gammafs}. \\
&  \leq   \frac{2}{\Delta t} C^2_{\Gamma}h^{2p+2} \| \phi \|^2_{p,\Gammafs}
+ \frac{g^2 \Delta t}{2}C^2_{\Gamma}h^{2p+2} \| \eta \|_{p,\Gammafs} \\
\end{alignat}
\end{subequations} 
\end{comment}
\begin{align}
 \vertiii{\bvphi - \bvphi^h}^2_m 
  \leq &
  \frac{2}{\Delta t} C^2_{\Gamma}h^{2p+2} \| \phi \|^2_{p,\Gammafs} \nn\\
&+ \frac{g^2 \Delta t}{2}C^2_{\Gamma}h^{2p+2} \| \eta \|_{p,\Gammafs}.
\end{align}

\end{document}

%% file: declarations.tex
\def\B#1{\mbox{\boldmath{$#1$}}}

\newcommand{\onequart}{\mbox{$\frac{1}{4}$}}

\newcommand{\nn}{\nonumber}

\newcommand{\bu}{\B{u}}

\newcommand{\bx}{\B{x}}
\newcommand{\by}{\B{y}}
\newcommand{\bn}{\B{n}}
\newcommand{\be}{\B{e}}

\newcommand{\bW}{\B{W}}

\newcommand{\bvphi}{\boldsymbol{\varPhi}}

\newcommand{\onehalf}{\mbox{$\frac{1}{2}$}}
\newcommand{\VV}{\mathcal{V}}
\newcommand{\WW}{\mathcal{W}}

\def\be{\begin{equation}}
\def\ee{\end{equation}}
\def\ba{\begin{array}}
\def\ea{\end{array}}
\def\bea{\begin{eqnarray}}
\def\eea{\end{eqnarray}}
\def\beas{\begin{eqnarray*}}
\def\eeas{\end{eqnarray*}}
\newcommand{\bseq}{\begin{subequations}}
\newcommand{\eseq}{\end{subequations}}

%% file: figs_bw/trace2_p1.tex
\begin{figure}[!htb]
\begin{center}
\begin{subfigure}[b]{0.39\textwidth}
\centering
\includegraphics[width=\textwidth]{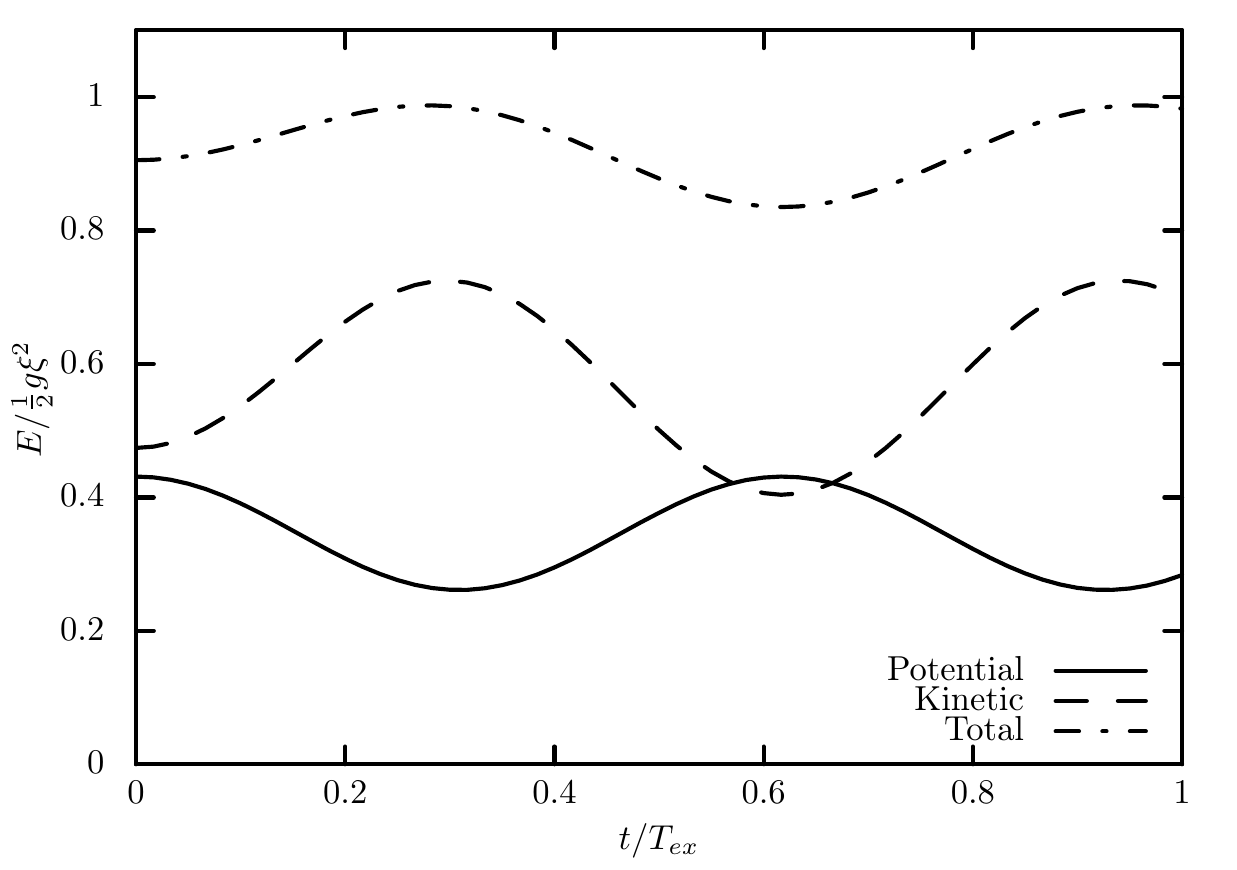}
\caption{6 x 6 mesh}
\end{subfigure}
\begin{subfigure}[b]{0.39\textwidth}
\centering
\includegraphics[width=\textwidth]{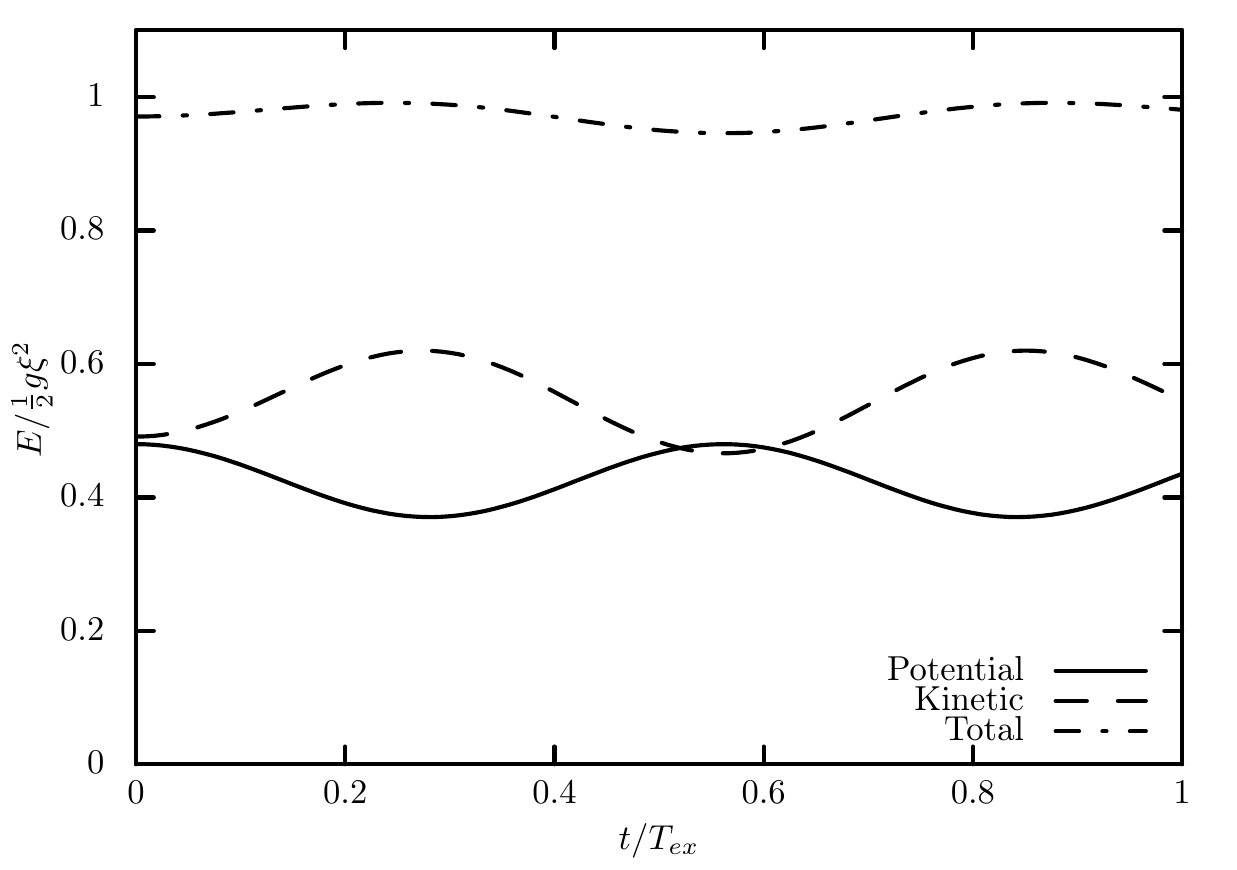}
\caption{12 x 12 mesh}
\end{subfigure}
\begin{subfigure}[b]{0.39\textwidth}
\centering
\includegraphics[width=\textwidth]{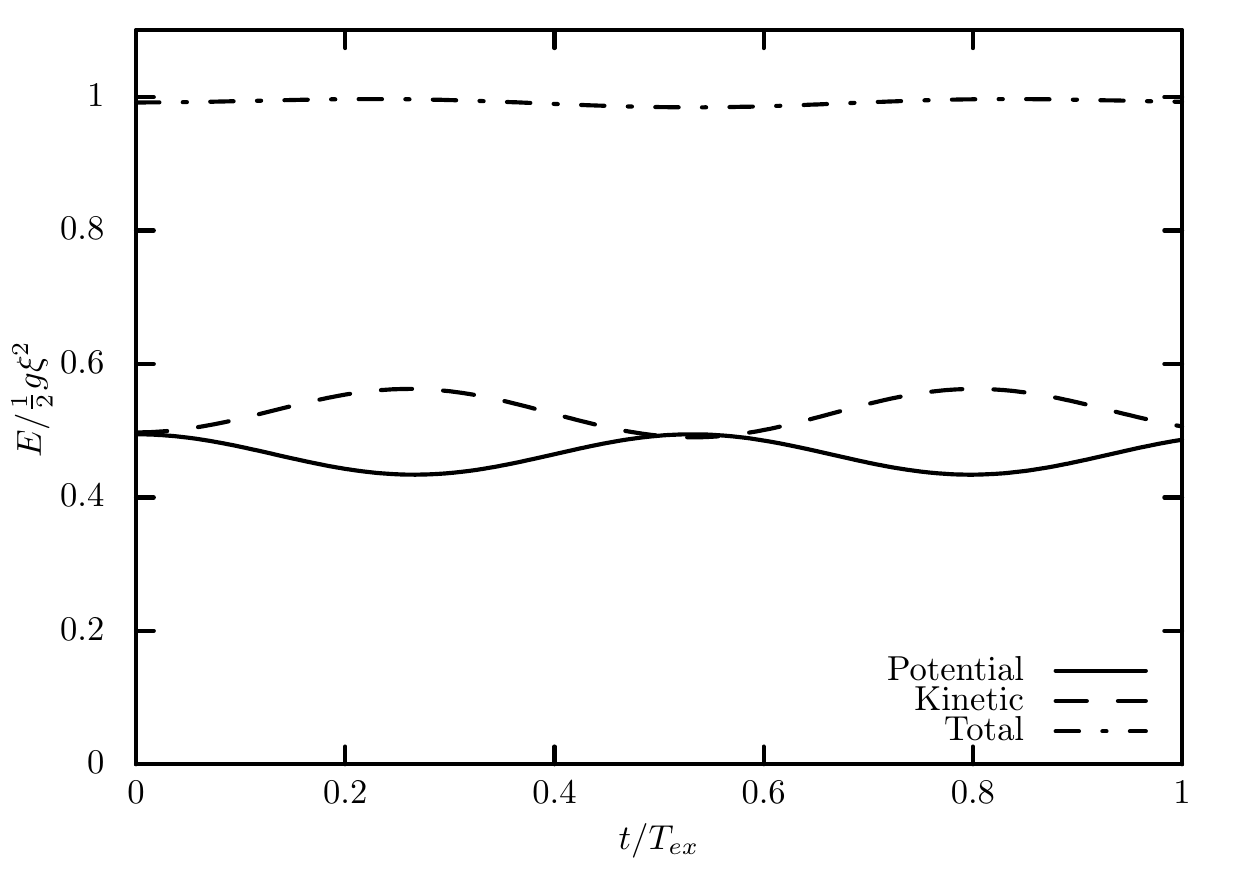}
\caption{24 x 24 mesh}
\end{subfigure}
\begin{subfigure}[b]{0.39\textwidth}
\centering
\includegraphics[width=\textwidth]{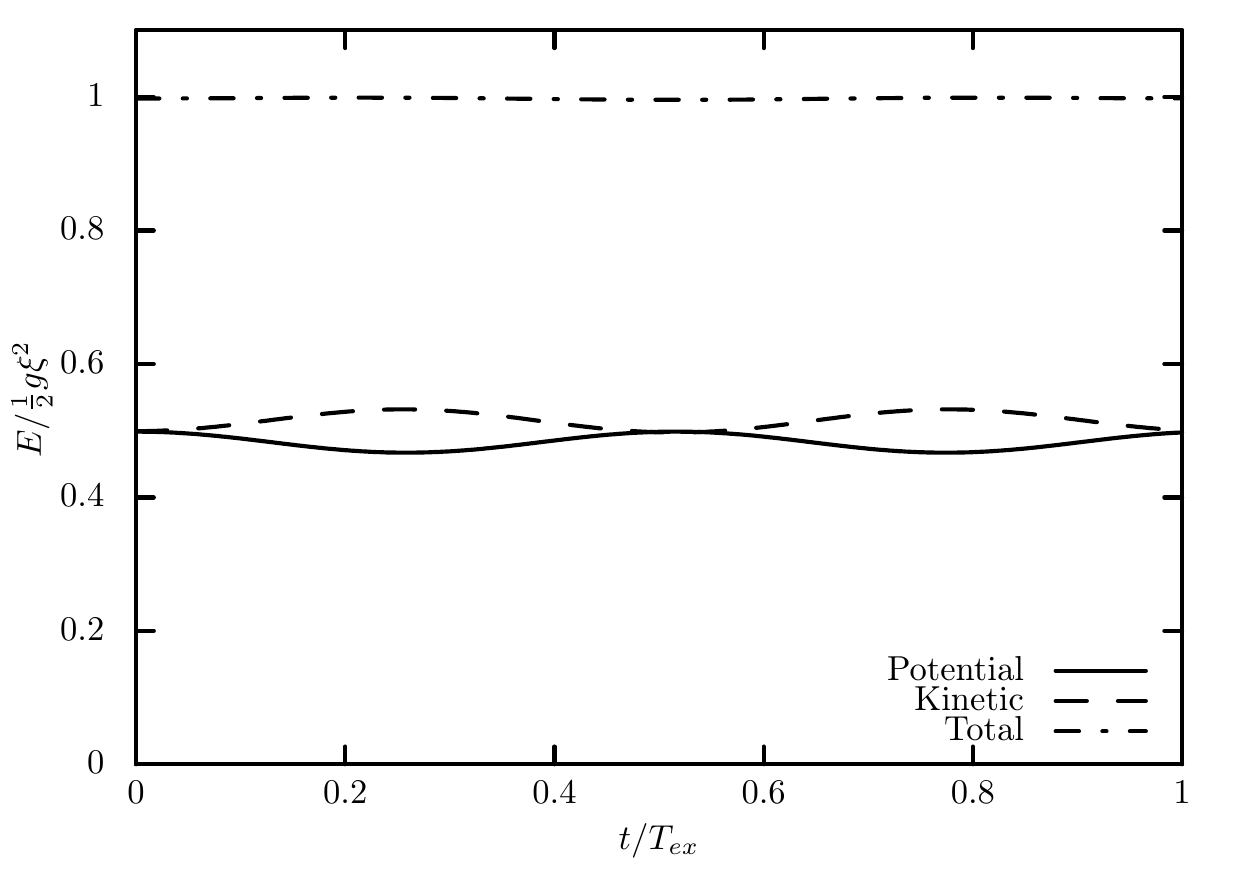}
\caption{48 x 48 mesh}
\end{subfigure}
\end{center}
\caption{Energy time trace for the segregated formulation with linear finite elements}
\label{fig:tracef3p1}
\end{figure}

\begin{figure}[!htb]
\begin{center}
\begin{subfigure}[b]{0.39\textwidth}
\centering
\includegraphics[width=\textwidth]{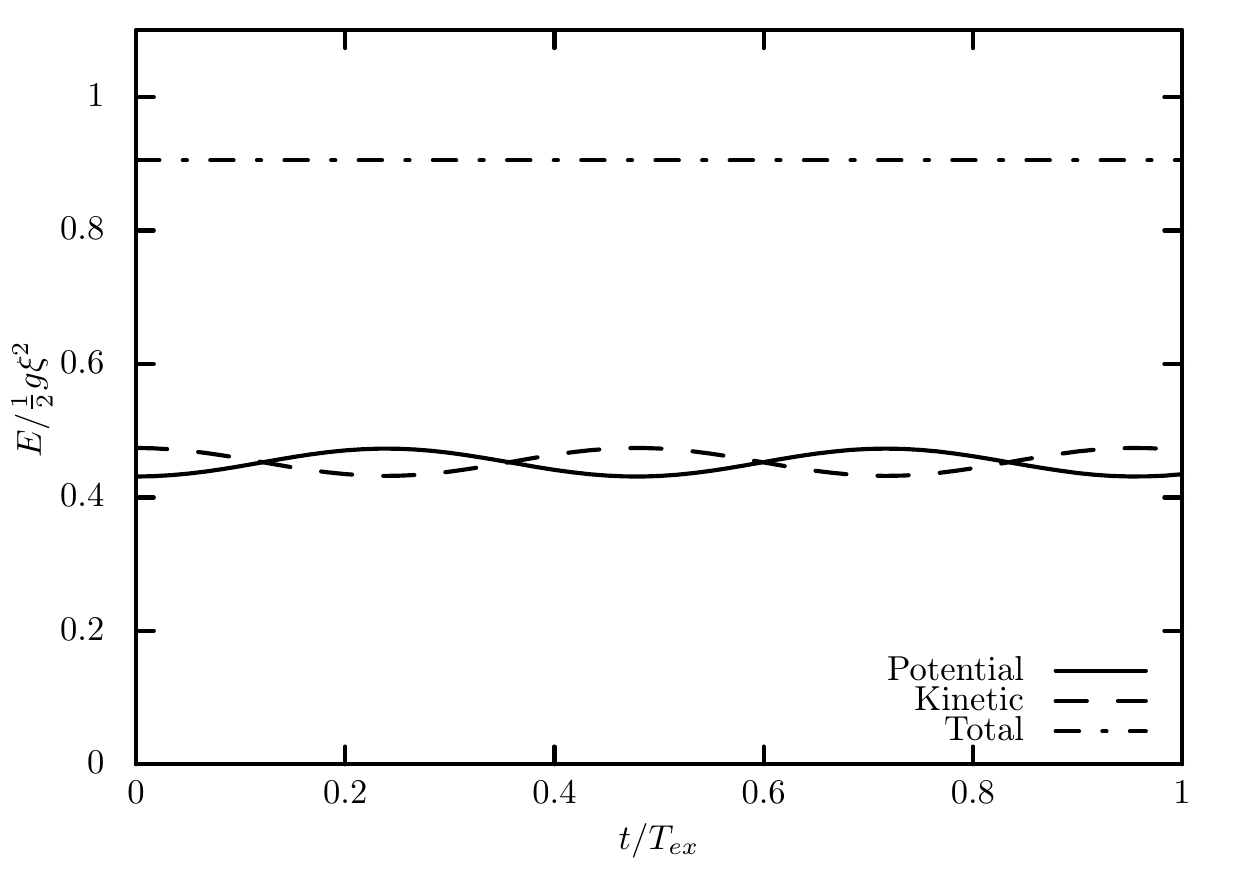}
\caption{6 x 6 mesh}
\end{subfigure}
\begin{subfigure}[b]{0.39\textwidth}
\centering
\includegraphics[width=\textwidth]{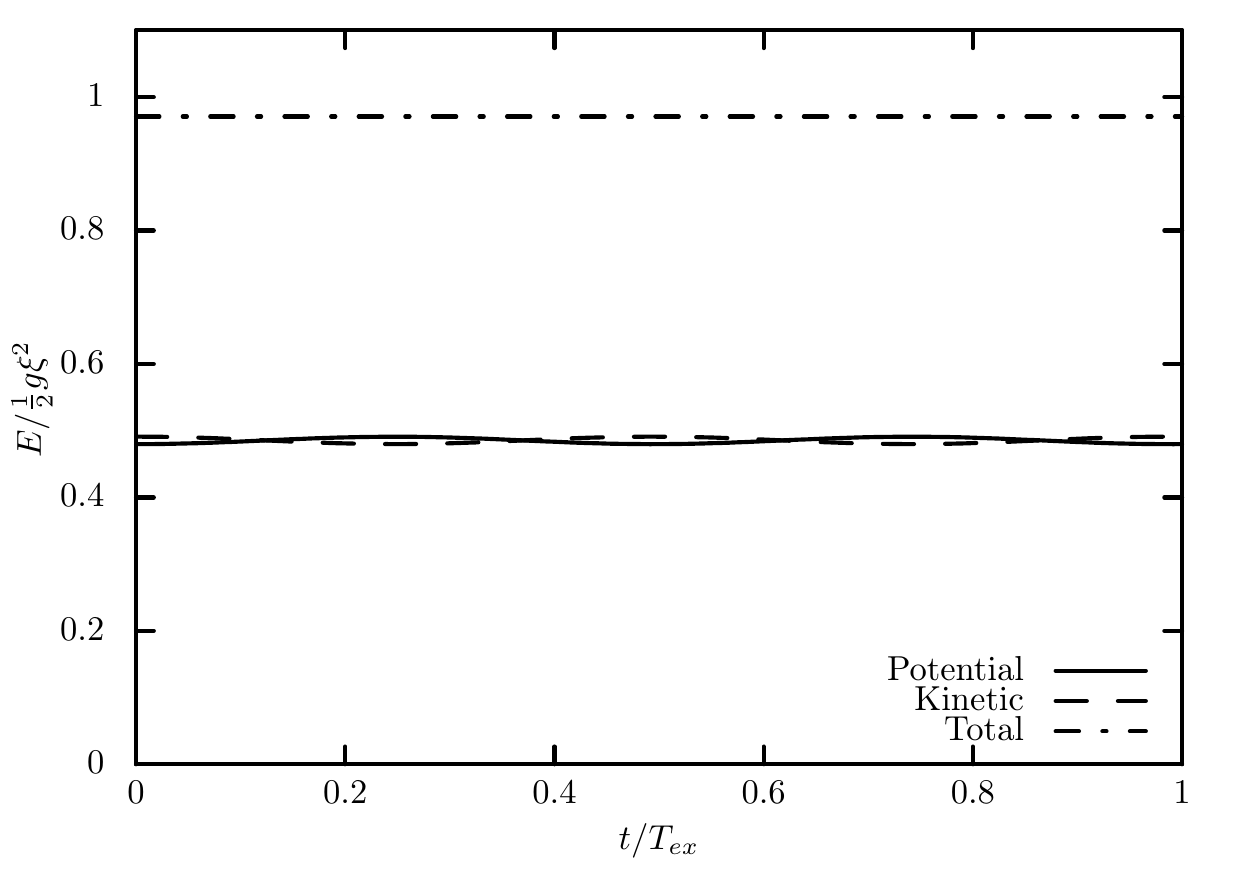}
\caption{12 x 12 mesh}
\end{subfigure}
\begin{subfigure}[b]{0.39\textwidth}
\centering
\includegraphics[width=\textwidth]{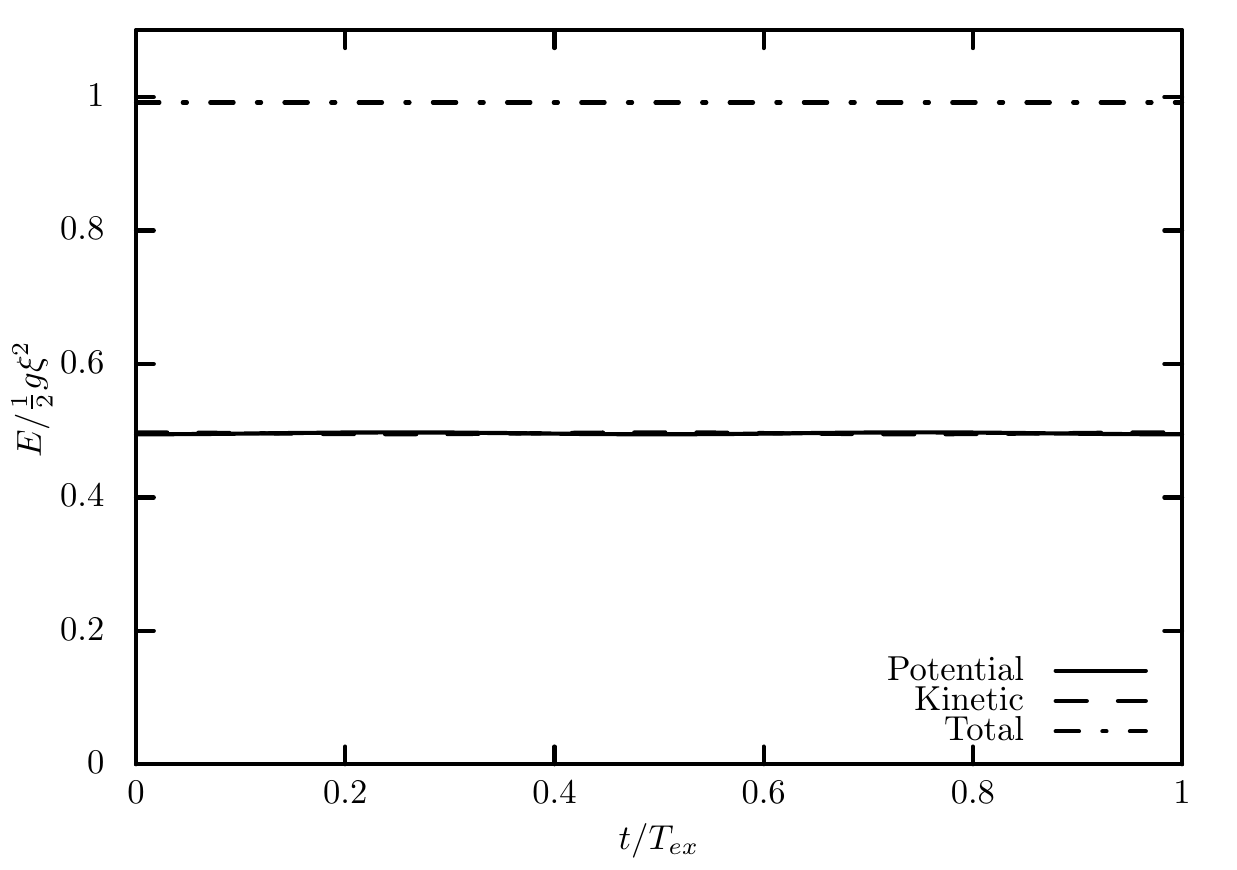}
\caption{24 x 24 mesh}
\end{subfigure}
\begin{subfigure}[b]{0.39\textwidth}
\centering
\includegraphics[width=\textwidth]{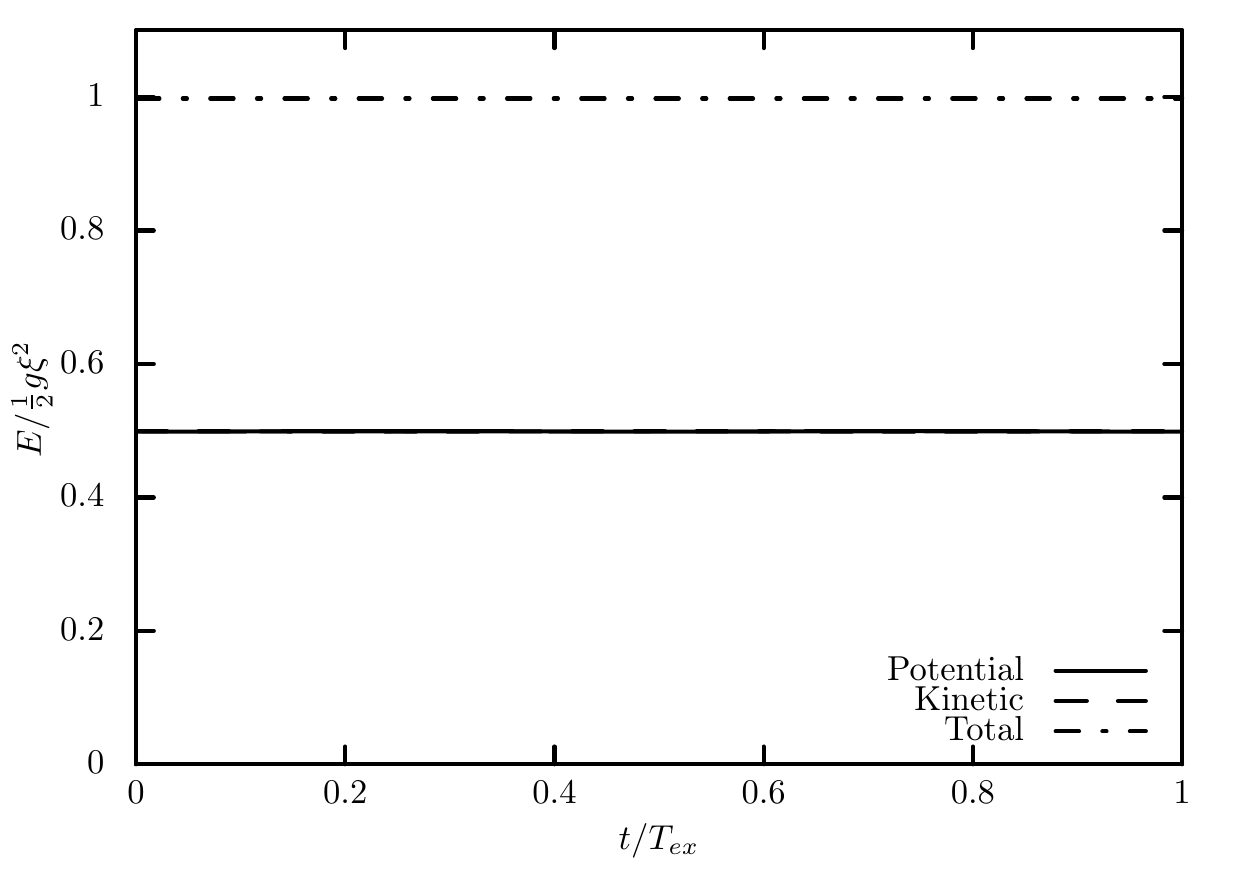}
\caption{48 x 48 mesh}
\end{subfigure}
\end{center}
\caption{Energy time trace for segregated formulation with linear finite elements}
\label{fig:tracef2p1}
\end{figure}

%% file: figs_bw/trace2_p2.tex
\begin{figure}[!htb]
\begin{center}
\begin{subfigure}[b]{0.39\textwidth}
\centering
\includegraphics[width=\textwidth]{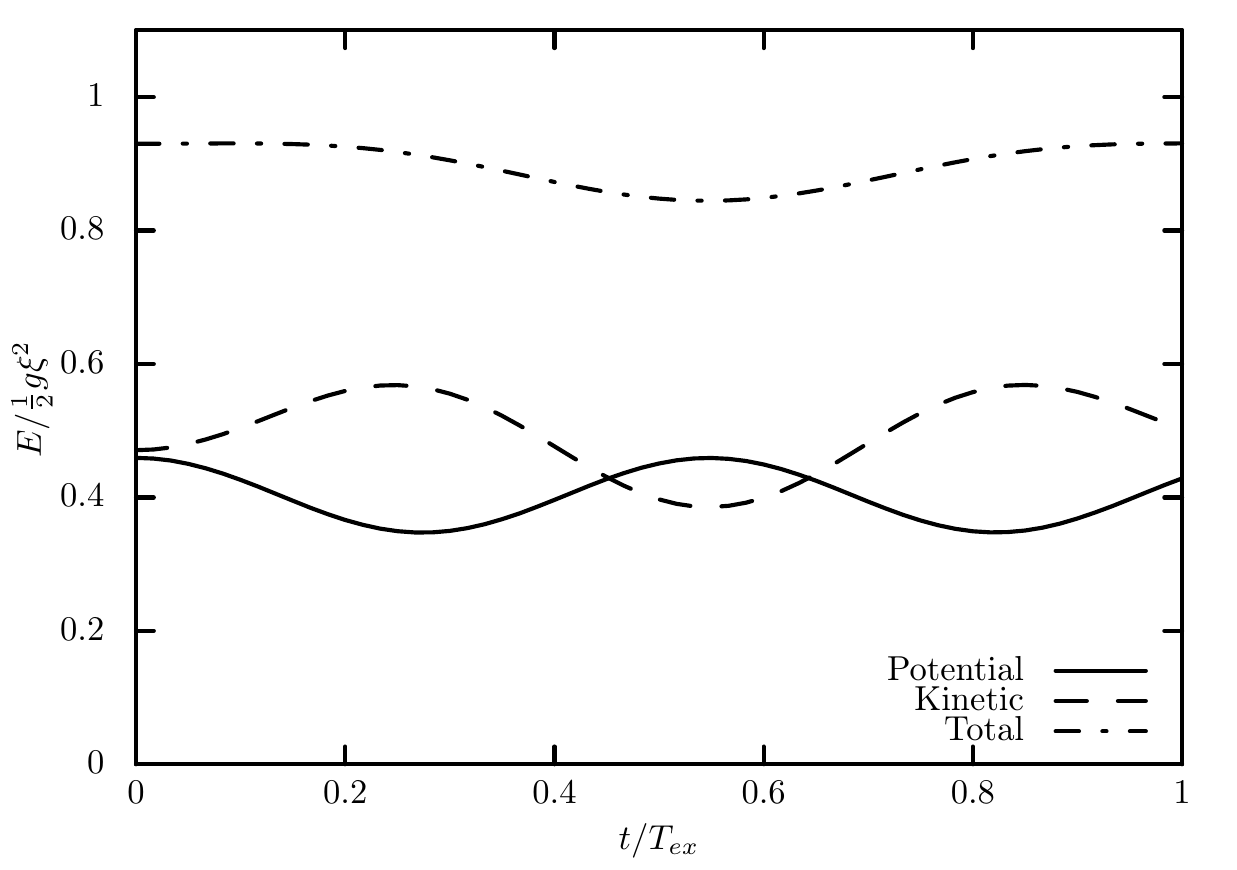}
\caption{3 x 3 mesh}
\end{subfigure}
\begin{subfigure}[b]{0.39\textwidth}
\centering
\includegraphics[width=\textwidth]{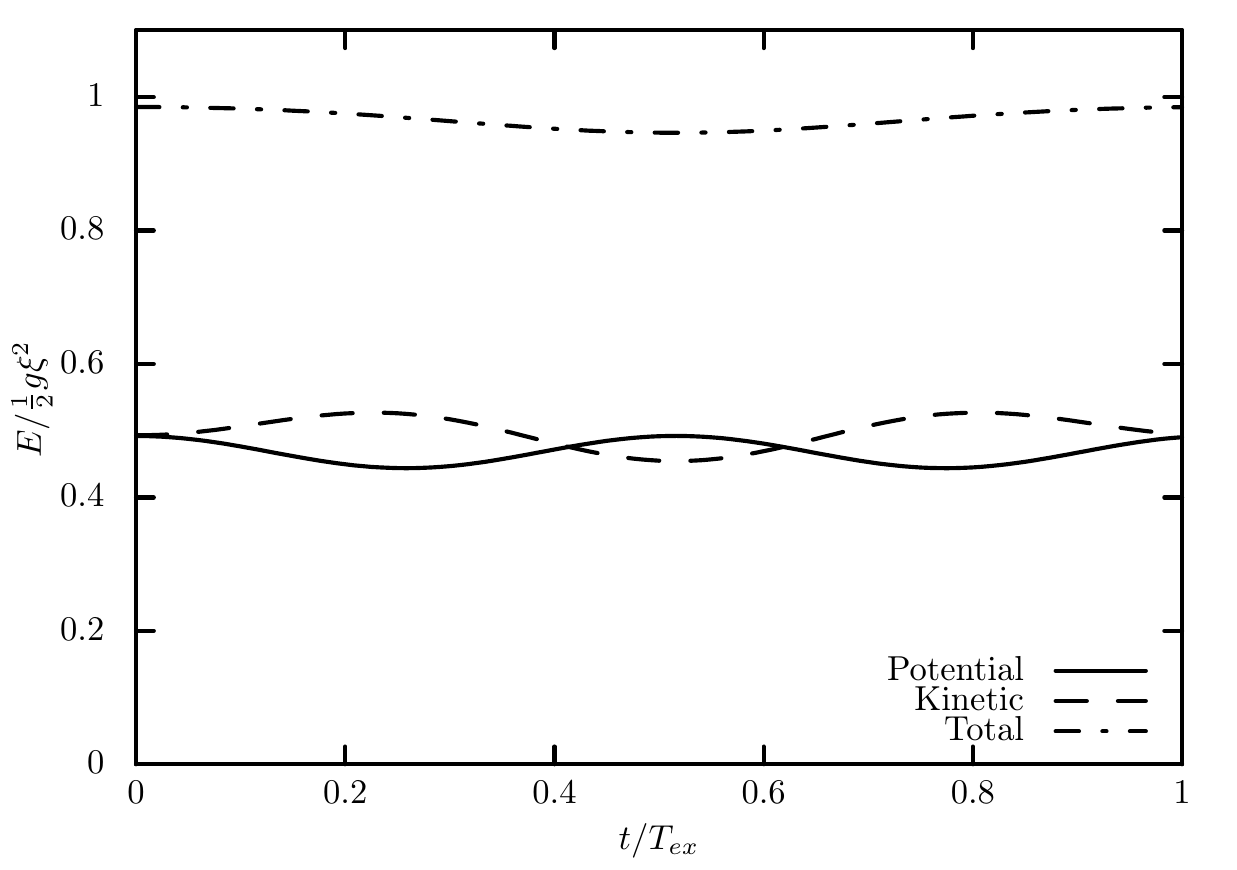}
\caption{6 x 6 mesh}
\end{subfigure}
\begin{subfigure}[b]{0.39\textwidth}
\centering
\includegraphics[width=\textwidth]{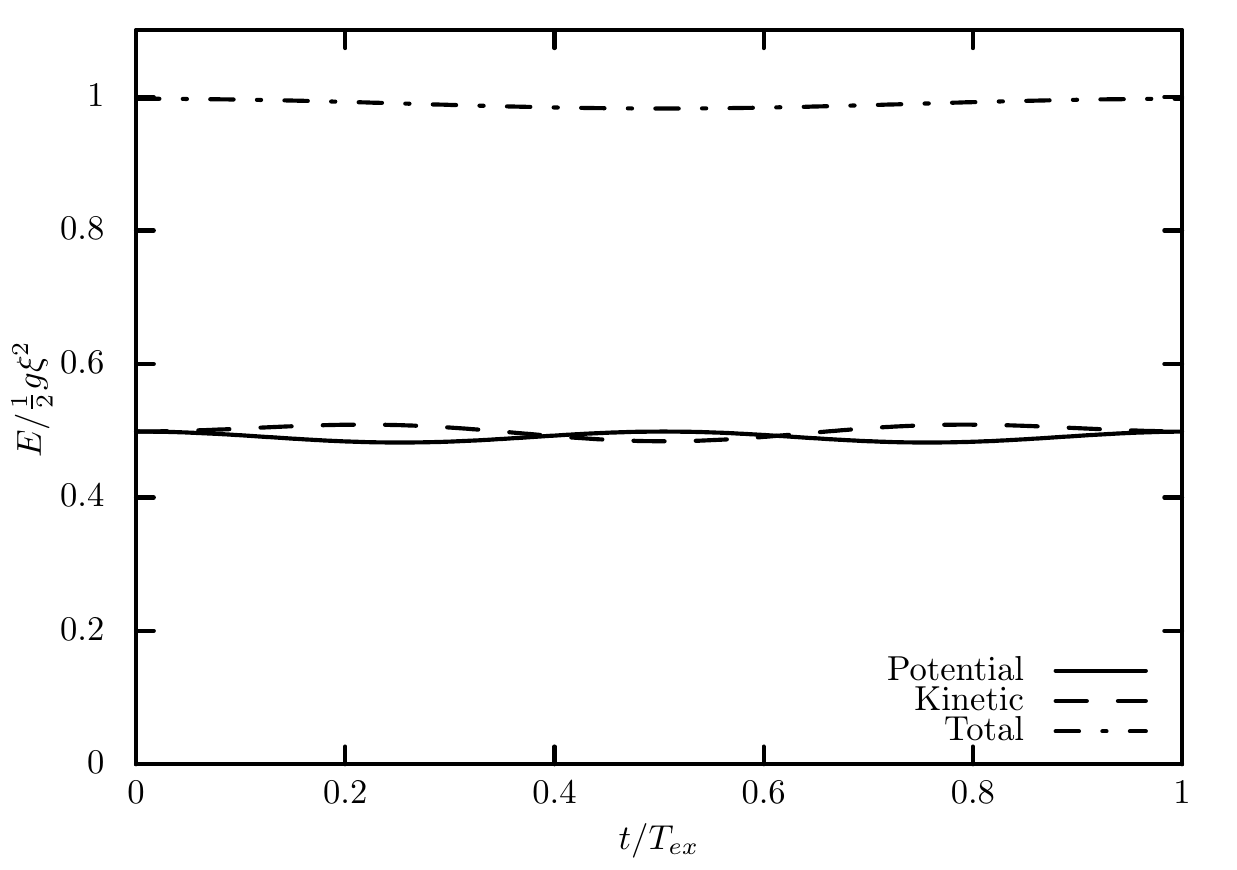}
\caption{12 x 12 mesh}
\end{subfigure}
\begin{subfigure}[b]{0.39\textwidth}
\centering
\includegraphics[width=\textwidth]{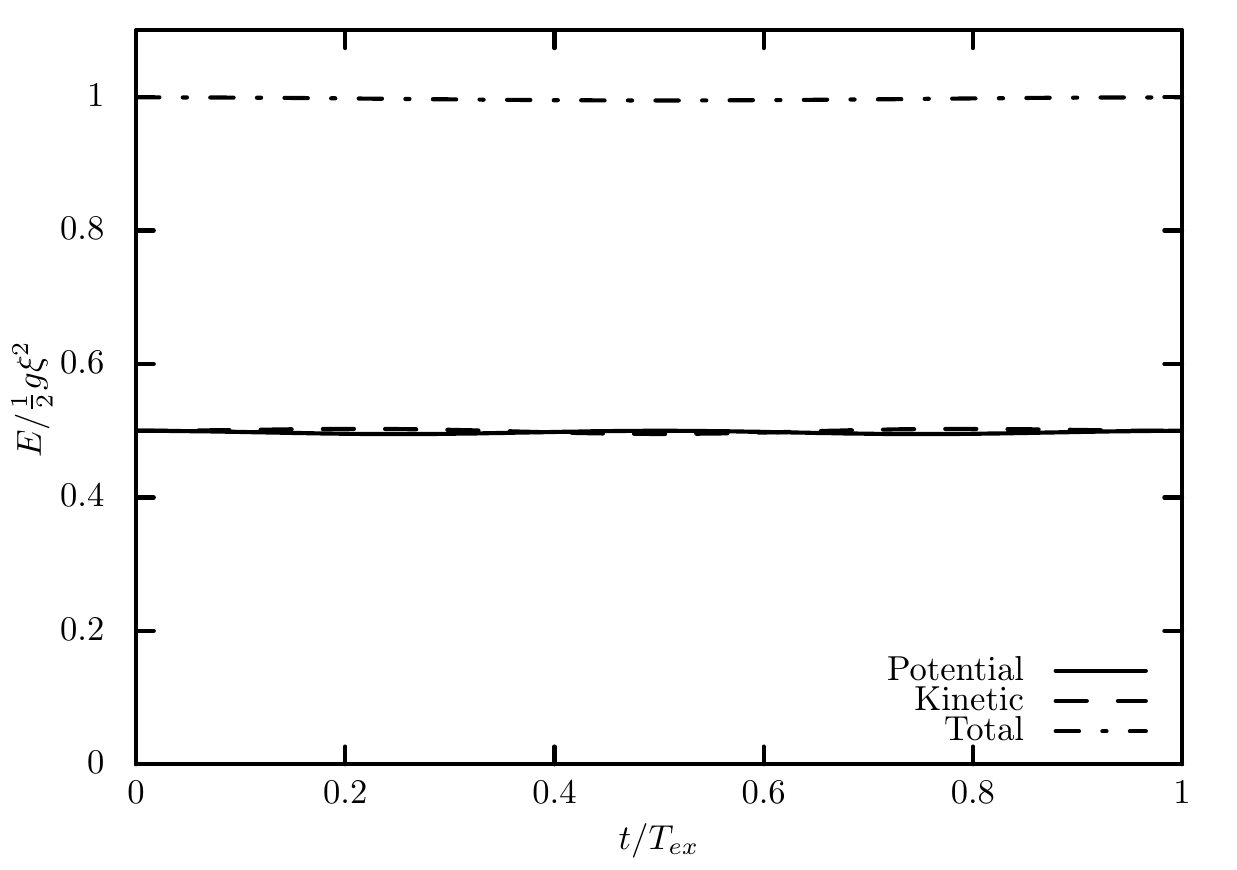}
\caption{24 x 24 mesh}
\end{subfigure}
\end{center}
\caption{Energy time trace for the segregated formulation with quadratic finite elements}
\label{fig:tracef3p2}
\end{figure}

\begin{figure}[!htb]
\begin{center}
\begin{subfigure}[b]{0.39\textwidth}
\centering
\includegraphics[width=\textwidth]{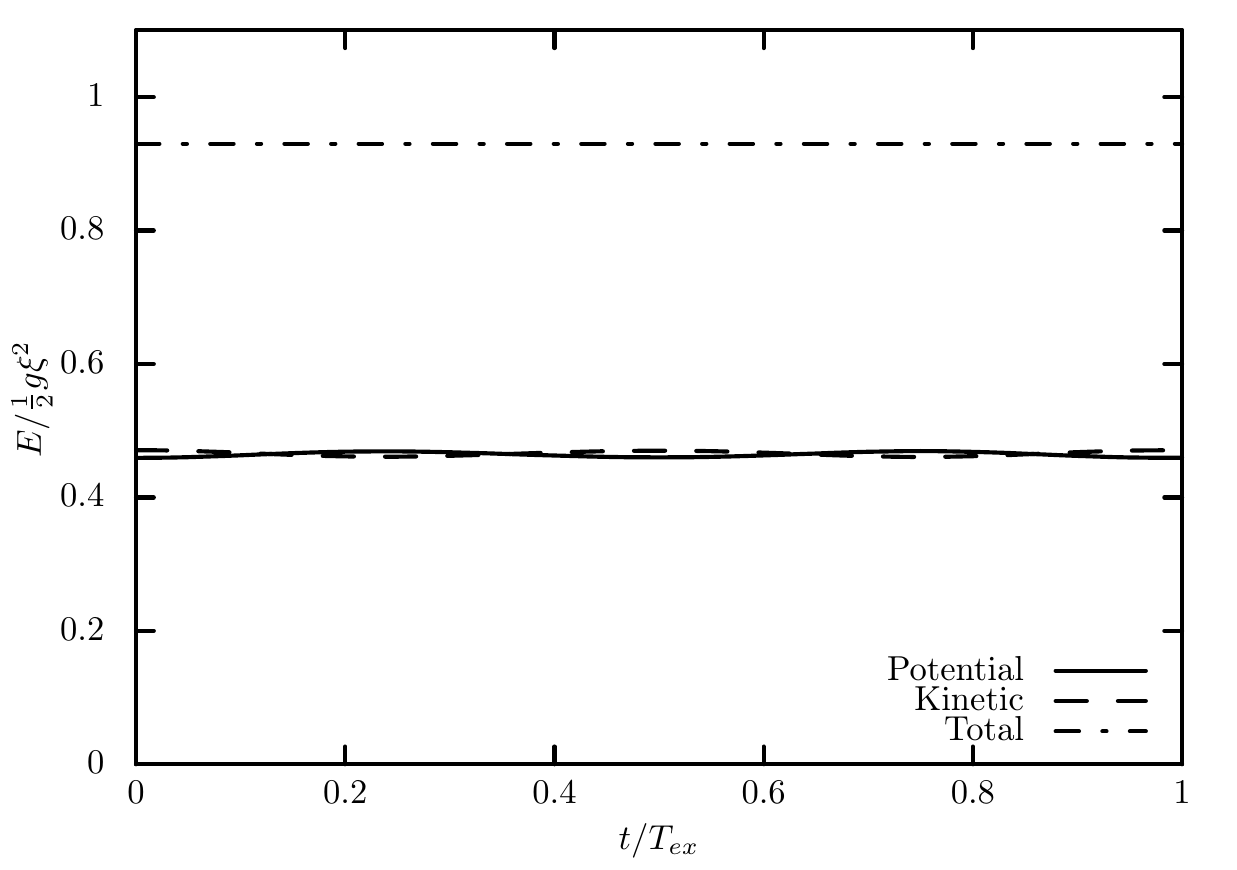}
\caption{3 x 3 mesh}
\end{subfigure}
\begin{subfigure}[b]{0.39\textwidth}
\centering
\includegraphics[width=\textwidth]{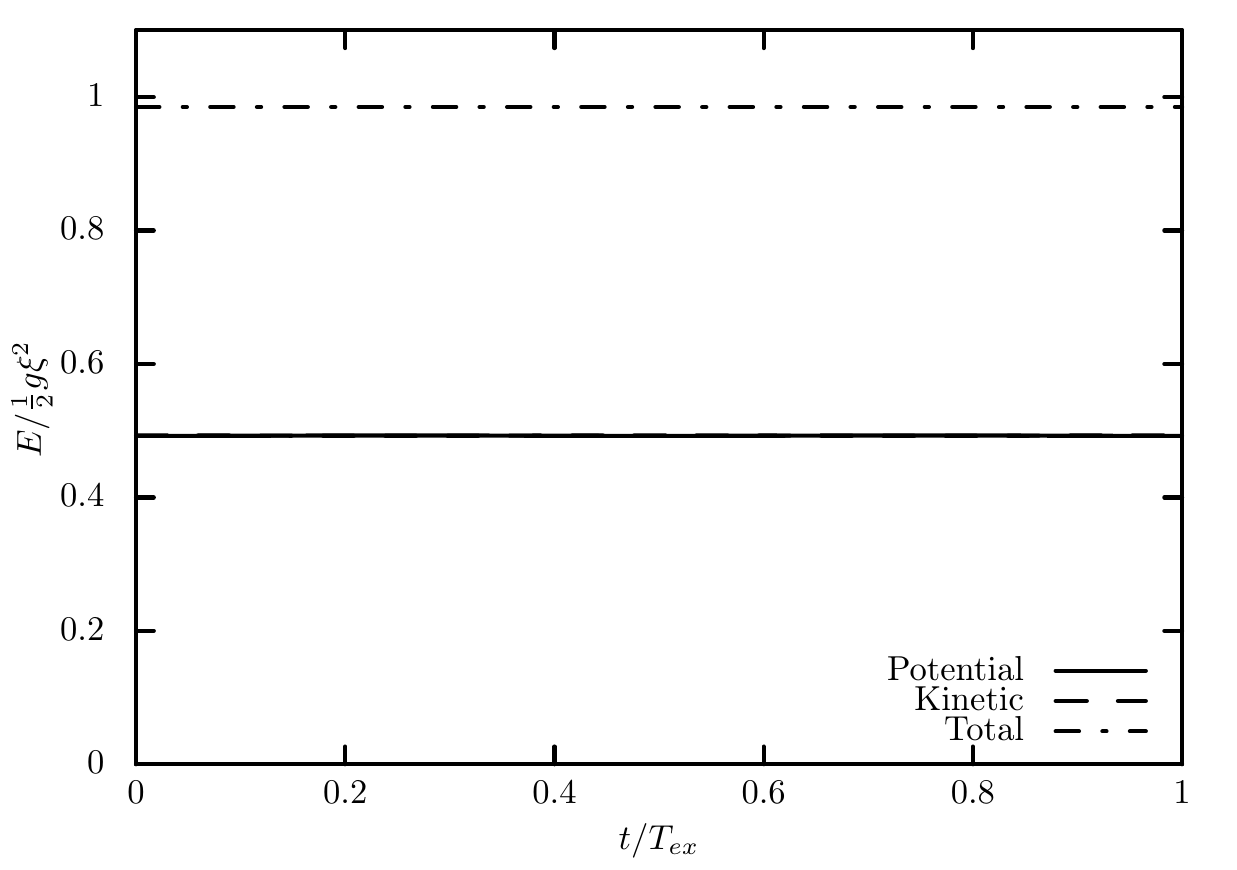}
\caption{6 x 6 mesh}
\end{subfigure}
\begin{subfigure}[b]{0.39\textwidth}
\centering
\includegraphics[width=\textwidth]{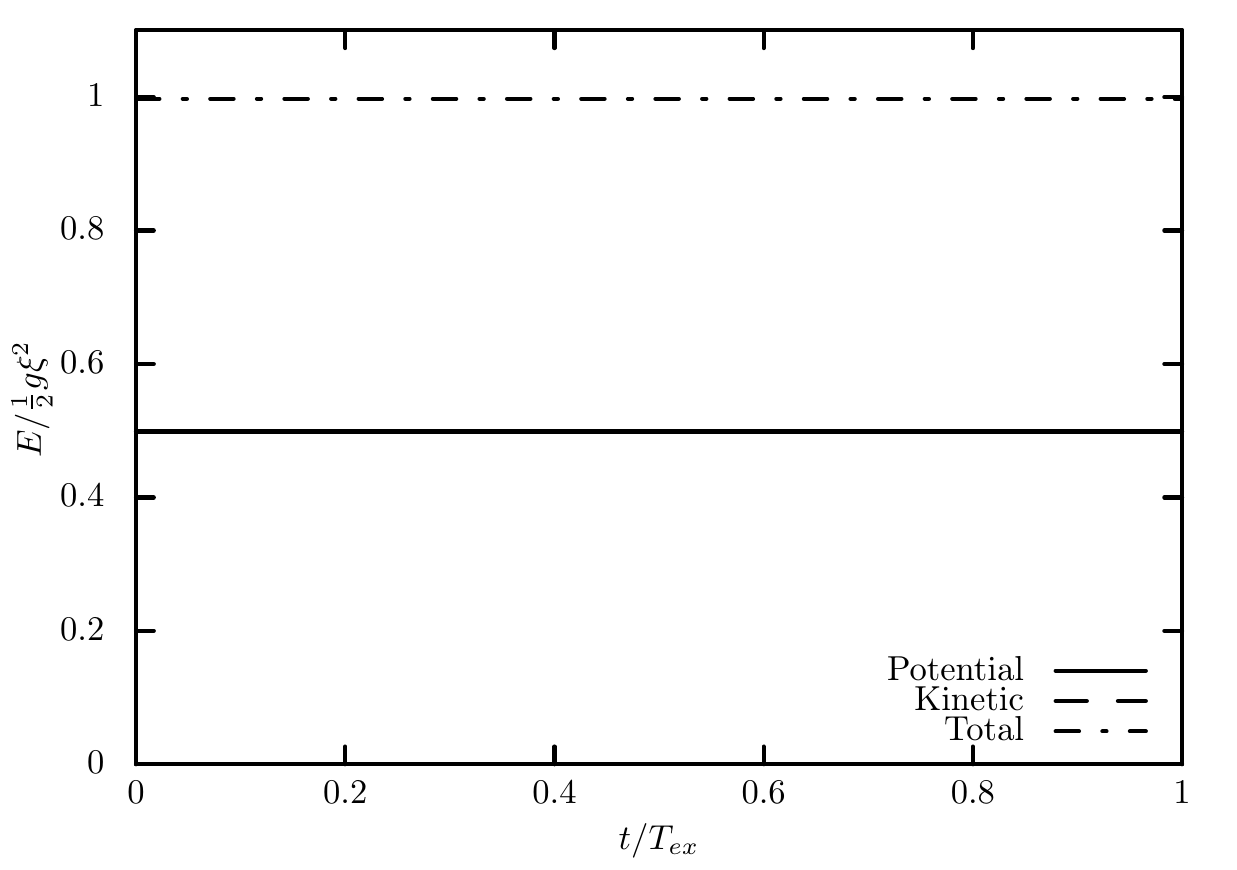}
\caption{12 x 12 mesh}
\end{subfigure}
\begin{subfigure}[b]{0.39\textwidth}
\centering
\includegraphics[width=\textwidth]{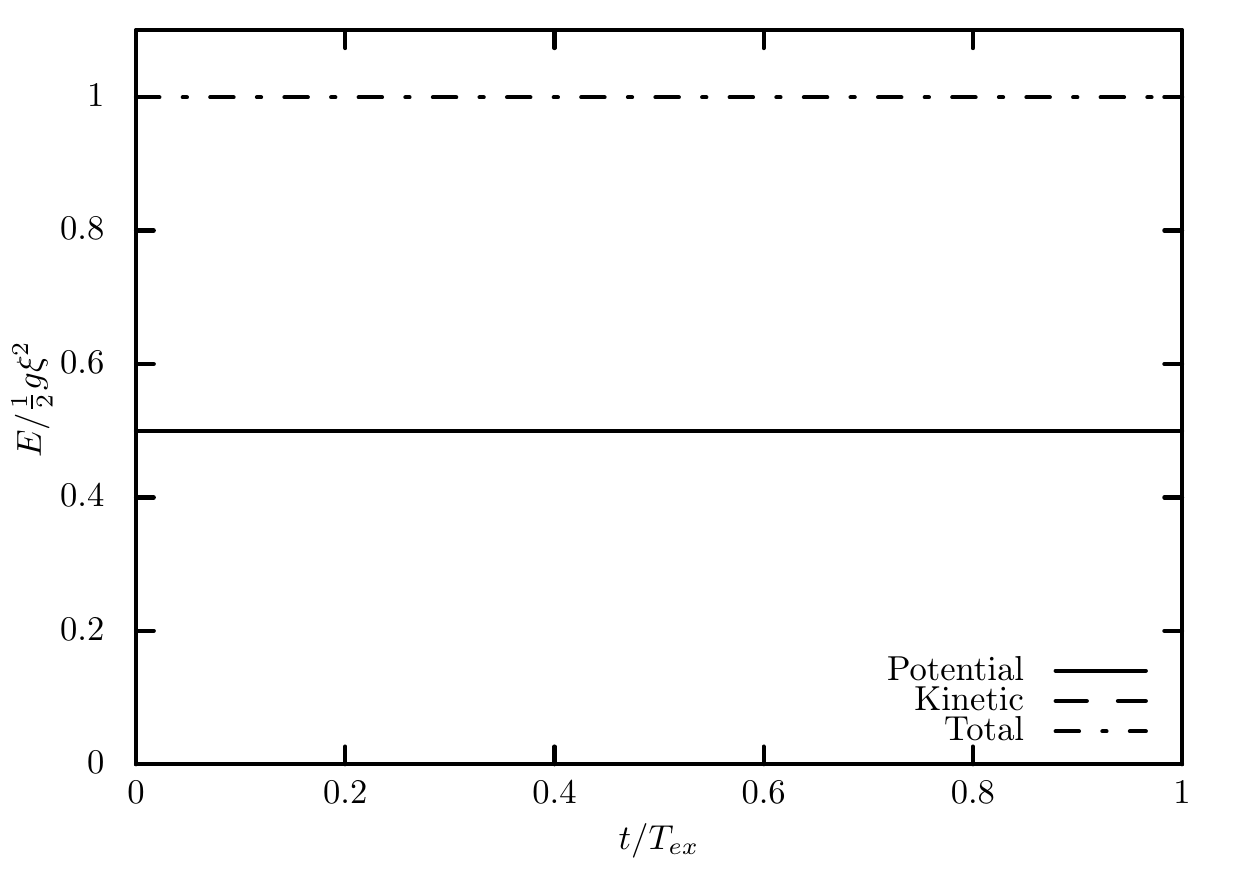}
\caption{24 x 24 mesh}
\end{subfigure}
\end{center}
\caption{Energy time trace for the monolithic formulation with quadratic finite elements}
\label{fig:tracef2p2}
\end{figure}

%% file: paper1.bbl
\begin{thebibliography}{10}

\bibitem{Pinkster79}
J.A. Pinkster.
\newblock {Mean and low-frequency wave drifting forces on floating structures}.
\newblock {\em {Ocean Engineering}}, {6}:{593--615}, {1979}.

\bibitem{Dawson77}
C.W. Dawson.
\newblock A practical computer method for solving ship-wave problems.
\newblock {\em Proc. 2nd Intl. Conf. on Numerical Ship Hydrodynamics}, pages
  30--38, 1977.

\bibitem{STF70}
N.~Salvesen, E.O. Tuck, and O.~Faltinsen.
\newblock Ship motions and sea loads.
\newblock {\em Trans. SNAME}, 78:250--287, 1970.

\bibitem{CRW93}
R.~Coifman, V.~Rokhlin, and S.~Wandzura.
\newblock The fast multipole method for the wave equation: A pedestrian
  prescription.
\newblock {\em IEEE Antennas and Propagation Magazine}, 35:7--12, 1993.

\bibitem{HuCoBa04}
T.J.R. Hughes, J.A. Cottrell, and Y.~Bazilevs.
\newblock Isogeometric analysis: {CAD}, finite elements, {NURBS}, exact
  geometry, and mesh refinement.
\newblock {\em Computer Methods in Applied Mechanics and Engineering},
  194:4135--4195, 2005.

\bibitem{CoHuBa09}
J.A. Cottrell, T.J.R. Hughes, and Y.~Bazilevs.
\newblock {\em Isogeometric Analysis: Toward Integration of CAD and FEA}.
\newblock Wiley, Chichester, 2009.

\bibitem{BACHH07}
I.~Akkerman, Y.~Bazilevs, V.M. Calo, T.J.R. Hughes, and S.~Hulshoff.
\newblock The role of continuity in residual-based variational multiscale
  modeling of turbulence.
\newblock {\em Computational {M}echanics}, 41:371--378, 2008.

\bibitem{CoHuRe07}
J.A. Cottrell, T.J.R. Hughes, and A.~Reali.
\newblock Studies of refinement and continuity in isogeometric structural
  analysis.
\newblock {\em Computer Methods in Applied Mechanics and Engineering},
  196:4160--4183, 2007.

\bibitem{EBBH09}
J.A. Evans, Y.~Bazilevs, I.~Babu{\v{s}}ka, and T.J.R. Hughes.
\newblock n-widths, sup--infs, and optimality ratios for the k-version of the
  isogeometric finite element method.
\newblock {\em Computer Methods in Applied Mechanics and Engineering},
  198:1726--1741, 2009.

\bibitem{WuTaylor94}
G.X. Wu and R.E. Taylor.
\newblock Finite element analysis 2-{D}imensional nonlinear transient
  water-waves.
\newblock {\em Applied Ocean Research}, {16}:{363--372}, {1994}.

\bibitem{kim1999finite}
J.W. Kim and K.J. Bai.
\newblock A finite element method for two-dimensional water-wave problems.
\newblock {\em International journal for numerical methods in fluids},
  30:105--122, 1999.

\bibitem{phdWesthuis}
J.-H. Westhuis.
\newblock {\em The numerical simulation of nonlinear waves in a hydrodynamic
  model test basin}.
\newblock PhD thesis, University Twente, 2001.

\bibitem{ZZ92}
O.C. Zienkiewicz and J.Z. Zhu.
\newblock The superconvergent patch recovery and a posteriori error estimates.
  part 1: The recovery technique.
\newblock {\em International Journal for Numerical Methods in Engineering},
  33:1331--1364, 1992.

\bibitem{WuMaTaylor98}
G.X. Wu, Q.W. Ma, and R.E. Taylor.
\newblock {Numerical simulation of sloshing waves in a 3D tank based on a
  finite element method}.
\newblock {\em Applied Ocean Research}, {20}:{337--355}, {1998}.

\bibitem{KHKB2005}
J.H. Kyoung, S.Y. Hong, J.W. Kim, and K.J. Bai.
\newblock {Finite-element computation of wave impact load due to a violent
  sloshing}.
\newblock {\em {Ocean Engineering}}, {32}:{2020--2039}, {2005}.

\bibitem{BCCK05}
K.J. Bai, S.M. Choo, S.K. Chung, and D.Y. Kim.
\newblock {Numerical solutions for nonlinear free surface flows by finite
  element methods}.
\newblock {\em Applied Mathematics and Computation}, {163}:{941--959}, {2005}.

\bibitem{KKEB2006}
J.W. Kim, J.H. Kyoung, R.C. Ertekin, and K.J. Bai.
\newblock {Finite-element computation of wave-structure interaction between
  steep Stokes waves and vertical cylinders}.
\newblock {\em Journal of waterway port coastal and ocean engineering-ASCE},
  {132}:{337--347}, {2006}.

\bibitem{EiAk17i}
M.F.P. ten Eikelder and I.~Akkerman.
\newblock {Correct energy evolution of stabilized formulations: The relation
  between VMS, SUPG and GLS via dynamic orthogonal small--scales and
  isogeometric analysis. I: The convective--diffusive context}.
\newblock {\em Computer Methods in Applied Mechanics and Engineering},
  331:259--280, 2018.

\bibitem{Evans13unsteadyNS}
J.A. Evans and T.J.R. Hughes.
\newblock Isogeometric divergence-conforming {B}-splines for the unsteady
  {N}avier--{S}tokes equations.
\newblock {\em Journal of Computational Physics}, 241:141--167, 2013.

\bibitem{EiAk17ii}
M.F.P. ten Eikelder and I.~Akkerman.
\newblock {Correct energy evolution of stabilized formulations: The relation
  between VMS, SUPG and GLS via dynamic orthogonal small--scales and
  isogeometric analysis. II: The incompressible Navier--Stokes equations}.
\newblock {\em Computer Methods in Applied Mechanics and Engineering},
  340:1135--1154, 2018.

\bibitem{AkEik19}
I.~Akkerman and M.F.P. ten Eikelder.
\newblock Toward free-surface flow simulations with correct energy evolution:
  An isogeometric level-set approach with monolithic time-integration.
\newblock {\em Computers \& Fluids}, 181:77 -- 89, 2019.

\bibitem{Tspline2010}
Y.~Bazilevs, V.M. Calo, J.A. Cottrell, J.A. Evans, T.J.R. Hughes, S.~Lipton,
  M.A. Scott, and T.W. Sederberg.
\newblock Isogeometric analysis using t-splines.
\newblock {\em Computer Methods in Applied Mechanics and Engineering},
  199:229--263, 2010.

\bibitem{LRspline2013}
T.~Dokken, T.~Lyche, and K.F. Pettersen.
\newblock Polynomial splines over locally refined box-partitions.
\newblock {\em Computer Aided Geometric Design}, 30:331--356, 2013.

\bibitem{USpline2019}
D.C. Thomas, L.~Engvall, S.~Schmidt, K.~Tew, and M.A. Scott.
\newblock U-splines: Splines over unstructured meshes.
\newblock coreform.com/usplines.
\newblock Accessed: 2019-02-06.

\bibitem{ErnGuermondFEM2004}
A.~Ern and J.-L. Guermond.
\newblock {\em {Theory and Practice of Finite Elements}}.
\newblock Appl. Math. Sci. 159, Springer-Verlag, New York, 2004.

\bibitem{Ciarlet78}
Ph.G. Ciarlet.
\newblock {\em The Finite Element Method for Elliptic Problems}.
\newblock North Holland, 1978.

\bibitem{BBVCHS06}
Y.~Bazilevs, L.~Beirao~da Veiga, J.~A. Cottrell, T.J.R. Hughes, and
  G.~Sangalli.
\newblock Isogeometric analysis: approximation, stability and error estimates
  for h-refined meshes.
\newblock {\em Mathematical Models and Methods in Applied Sciences},
  16:1031--1090, 2006.

\bibitem{mfem-library}
{MFEM}: Modular finite element methods library.
\newblock mfem.org.

\bibitem{VisIt}
H.~Childs {\it et al}.
\newblock {VisIt: An End-User Tool For Visualizing and Analyzing Very Large
  Data}.
\newblock In {\em {High Performance Visualization--Enabling Extreme-Scale
  Scientific Insight}}, pages 357--372. Oct 2012.

\end{thebibliography}
